\documentclass[12pt]{article}
\usepackage{a4}
\usepackage{amsmath,amsthm,amssymb}
\usepackage[latin3]{inputenc}
\usepackage{epsfig}
\usepackage{fancyhdr}
 \oddsidemargin=\evensidemargin
\newcommand{\bitem}{\begin{itemize}}
\newcommand{\eitem}{\end{itemize}}
\newcommand{\benum}{\begin{enumerate}}
\newcommand{\eenum}{\end{enumerate}}

\newcommand{\Plim}{{\rm P}\!\lim\limits_{\!\!\!\!n\ri \infty} }
\newcommand{\Plimk}{{\rm P}\!\lim\limits_{\!\!\!\!k\ri \infty} }
\newcommand{\dimo}{{\it\bf Proof}.}
\newtheorem{theorem}{Theorem}[section]
\newtheorem{example}[theorem]{Example}
\newtheorem{remark}[theorem]{Remark}
\newtheorem{corollary}[theorem]{Corollary}
\newtheorem{lemma}[theorem]{Lemma}
\newtheorem{proposition}[theorem]{Proposition}
\newtheorem{definition}[theorem]{Definition}
\newcommand{\bex}{\begin{example}}
\newcommand{\eex}{\end{example}}
\newcommand{\beq}{\begin{equation} }
\newcommand{\eeq}{\end{equation} }
\newcommand{\n}{\noindent}

\newcommand{\N}{{\rm{I\! N}}}
\newcommand{\R}{{\rm{I\! R}}}
\newcommand{\ri}{\rightarrow}
\newcommand{\DN}{\Delta_{j} N^{(q)}}
\newcommand{\DNa}{\Delta_{j} N^{(1)}}
\newcommand{\DNb}{\Delta_{j} N^{(2)}}
\newcommand{\DW}{\Delta_{j} W^{(q)}}
\newcommand{\DXa}{\Delta_{j} X^{(1)}}
\newcommand{\DXb}{\Delta_{j} X^{(2)}}

\newcommand{\DDa}{\Delta_{j} D^{(1)}}
\newcommand{\DDb}{\Delta_{j} D^{(2)}}
\newcommand{\Da}{D^{(1)}}
\newcommand{\Db}{D^{(2)}}
\newcommand{\D}{D^{(q)}}
\newcommand{\DDq}{(\Delta_{j} D^{(q)})^2}

\newcommand{\DX}{\Delta_{j} X^{(q)}}
\newcommand{\DXq}{(\Delta_{j} X^{(q)})^2}
\newcommand{\DXdq}{(\Delta_{j} \tilde J_2^{(q)})^2}
\newcommand{\DXd}{(\Delta_{j} \tilde J_2^{(q)})^2}
\newcommand{\DY}{\Delta_{j} Y^{(q)}}
\newcommand{\DD}{\Delta_{j} D^{(q)}}

\newcommand{\DJu}{\Delta_{j} J_1^{(q)}}
\newcommand{\rh}{r_h}
\newcommand{\mc}{2h\log \frac 1 h}
\newcommand{\sqrh}{\sqrt{r_h}}
\newcommand{\intIi}{\int_{t_{j-1}}^{t_j}}
\newcommand{\Ii}{]t_{j-1},t_j]}

\newcommand{\vm}{\vspace{0.5cm}}

\title{IDENTIFYING THE COVARIATION BETWEEN THE DIFFUSION PARTS AND THE CO-JUMPS GIVEN DISCRETE OBSERVATIONS}
\author{Fabio Gobbi\footnote{Dipartimento di Matematica  per le Decisioni,
Universit\`a degli Studi di Firenze}\quad and\quad Cecilia
Mancini\footnote{Dipartimento di Matematica per le Decisioni,
Universit\`a degli Studi di Firenze}}

\begin{document}

\maketitle

\begin{abstract}
In this paper we consider two semimartingales driven by diffusions
and jumps. We allow both for finite activity and for infinite
activity jump components. Given discrete observations we
disentangle the {\it integrated covariation} (the covariation
between the two diffusion parts, indicated by $IC$) from the
co-jumps. This has important applications to multiple assets price
modeling for forecasting, option pricing,
risk and credit risk management.\\
 An approach commonly used to estimate $IC$ is to take the sum of the cross
products of the two processes increments; however this estimator
can be highly biased in the presence of jumps, since it approaches
the quadratic covariation, which contains also the co-jumps. Our
estimator of $IC$ is based on a threshold (or truncation)
technique allowing to isolate all the jumps in the  finite
activity case and the jumps over the threshold
 in the  infinite activity case. We prove that the estimator is consistent in both cases as the number of observations
 increases to infinity.
Further, in presence of only finite activity jumps 1) $\hat{IC}$ is also asymptotically Gaussian;
 2) a joint CLT for $\hat{IC}$ and  threshold estimators of the integrated variances of
the single processes allows to reach consistent and asymptotically
Gaussian estimators even of the $\beta$s and of the correlation
coefficient among the diffusion parts of the two processes,
allowing a better measurement of their dependence; 3) thresholding
gives an estimate of $IC$ which is robust to the
asynchronicity of the observations. \\
We conduct a simulation study to check that the application of our technique is in fact informative
for values of  the step between the observations  large enough to avoid the typical
problems arising in presence of microstructure noises in the data,
and to asses the choice of the threshold parameters.\\

\end{abstract}

\section{Introduction}
We consider two state variables evolving as follows
$$
dX^{(q)}_{t}=
a^{(q)}_{t}dt+\sigma^{(q)}_{t}dW^{(q)}_{t}+dJ_{t}^{(q)}, \quad q=1,2
$$
for $t\in [0,T]$, $T<\infty$ fixed, where $a$ and $\sigma$ are
cadlag stochastic processes; $W^{(2)}_{t}=\rho_{t}
W^{(1)}_{t}+\sqrt{1-\rho_{t}^{2}}W^{(3)}_{t}$; $W^{(1)}$ and
$W^{(3)}$ are independent standard Brownian motions; and $J^{(1)}$
and $J^{(2)}$ are possibly correlated pure jump semimartingales.
Given discrete observations $X^{(1)}_{t_{j}}, X^{(2)}_{\nu_{i}},$
with observation times spanned on
$[0,T]$, 
we are interested in the separate identification of the {\it
integrated covariation}
$IC_T:=\int_0^T\rho_{t}\sigma^{(1)}_t\sigma^{(2)}_t dt$, between
the two diffusion parts, and of the {\it co-jumps}  $\Delta
J^{(1)}_{t}\Delta J^{(2)}_{t}$, the simultaneous jumps of
$X^{(1)}$ and $X^{(2)}$, where, for each $q=1,2$, $\Delta
J^{(q)}_{t}$ denotes the size $J^{(q)}_{t}- J^{(q)}_{t-}$ of the
jump occurred at time $t$.

The recent empirical interest on co-jumps in financial
econometrics 
is motivated by the problem of a
correct assets price model selection. This has important
consequences in forecasting, in option pricing, 
in portfolio risk management, and even in the credit risk
management, since a default of a firm is interpretable as a jump
in the firm value and contemporaneous defaults give a co-jump,
implying default dependence (contagion, \cite{ELV05}).

A commonly used approach to estimate
$\int_0^T\rho_{t}\sigma^{(1)}_t\sigma^{(2)}_t dt$ is to take synchronous and evenly-spaced observations
$X^{(1)}_{t_0},X^{(1)}_{t_1},...X^{(1)}_{t_n}, X^{(2)}_{t_0},X^{(2)}_{t_1},...X^{(2)}_{t_n}$, with $t_n=T,$ and
to consider  the sum
of cross products $\sum_{j=1}^{n}\Delta_j X^{(1)}\Delta_j X^{(2)}$,
where $\Delta_j X^{(q)}:= X^{(q)}_{t_{j}}-X^{(q)}_{t_{j-1}}$;
however this estimate can be highly biased when the processes
$X^{(q)}$ contain jumps; in fact, as $n\ri \infty,$ such a sum
approaches the global quadratic covariation
$$[X^{(1)},X^{(2)}]_{T}=\int_{0}^{T}\rho_{t}
\sigma_{t}^{(1)}\sigma_{t}^{(2)}dt+\sum_{0\leq t\leq T}\Delta
J^{(1)}_{t}\Delta J^{(2)}_{t},$$  which contains also the co-jumps.
To our aim it is crucial to single out the time intervals where the
jumps occurred.

A jump process $J$ is said to have {\it finite activity} (FA) when
a.s. only a finite number of jumps can occur in each finite time
interval. On the contrary $J$ is said to have {\it infinite
activity} (IA). In the special case where $J$ is Lévy and has IA
then {\it a.s.} infinitely many jumps occur in {\it each} finite
time interval.

Our estimator of $IC_T$ is based on a threshold criterion 
allowing to identify all the time intervals $]t_{j-1}, t_j]$ where
the path of a univariate semimartingale jumped, if the jump
component $J$ has FA, and the intervals where jumps over the
threshold occurred, if the discretely observed realization of $J$ has infinite activity.
Extending the application of the criterion to a bivariate framework allows to derive
 an asymptotically unbiased estimator of $IC_T$
as well as of the  co-jumps occurred up to time $T$. More
precisely
 we construct the following estimator
$$
\hat{IC}_{T,n}:= \sum_{j=1}^{n}
\Delta_{j}X^{(1)}1_{\{(\Delta_{j}X^{(1)})^{2}\leq
r(h)\}}\Delta_{j}X^{(2)}1_{\{(\Delta_{j}X^{(2)})^{2}\leq r(h)\}},
\quad h=T/n,
$$
 where only the variations under a given threshold function $r(h)$ are taken into account.
 The first main result of our paper
 is showing  the consistency to $IC_T$, as the number $n$ of observations tends to
 infinity.
 Not equally spaced but synchronous observations are allowed for such
 result. The second group of results is given in presence of only FA
 jumps. If we dispose of non-synchronous data we still reach consistency by modifying our estimator
 in a similar way of \cite{HayYos05} and \cite{HayKus08}. When observations are evenly spaced, we prove a joint CLT
delivering: 1. that $\hat{IC}_{T,n}$ is asymptotically Gaussian
and converes with speed
 $\sqrt{h}$, which extends results in
  \cite{BarShe04} who estimated $IC_T$ in
absence of jumps; 2. consistent and asymptotically Gaussian
estimators of the regression coefficients $\beta$s and of the
correlation coefficient
of the two continuous parts of processes $X^{(q)}$.\\
 In a further paper (\cite{GobMan07}) we
explore the speed of convergence of the estimator $\hat{IC}_{T,n}$
proposed here in the presence of, possibly correlated, infinite
activity Lévy jump processes $J^{(1)}$ and $J^{(2)}$, with
dependence structure described by a Lévy copula.

The threshold criterion originated in \cite{Man01} to separate the
diffusion and the jump parts of a univariate parametric
Poisson-Gaussian model. The criterion was shown to work even in
nonparametric frameworks in \cite{Man04},  \cite{mancini07} and
\cite{Jac08}. Potentially the threshold technique can be useful
in each context where disentangling the quadratic variation of a
signal has some importance to capture the contribution given by
the diffusive component of the model and the one given  by the
jump component.

The literature on non parametric inference for stochastic processes driven by
diffusions plus jumps, based on discrete
observations, is mainly devoted to univariate cases.
As for bivariate processes Barndorff-Nielsen and  Shephard
(\cite{BarShe04BPC}) and Jacod and Torodov (\cite{JacTor07}) explore
tests for the presence of co-jumps based on estimators
constructed basically using cross multipower variations. \\
We adopt the threshold method here since, at least in the finite
activity case, it is a more effective way to identify
(asymptotically) the intervals between consecutive observations where 
 jumps occurred.
In fact already in the univariate case the threshold estimator of
$IV^{(1)}$ is efficient (in the Cramer-Rao inequality lower bound
sense), the asymptotic standard estimation error being
$\sqrt{2\times IQ^{(1)}}$, where
$IQ^{(1)}:=\int_0^T(\sigma_{t}^{(1)})^4 dt$  \ (see
\cite{mancini07}), while the multipower variation estimators are
not efficient (the standard errors are all higher, see \cite{BarShe04PVandBPV},
\cite{Woe06} and the discussion in \cite{mancini07}, section 3.3,
the infimum being $\sqrt{2.609\times IQ^{(1)}}$).
 For the
bipower covariation  based estimator (BPC) of $IC_T$,
 a CLT  has been shown to hold only in absence
of jumps (\cite{BarShe04BPC}), in which case the standard error is a
function of $\rho_t, \sigma^{(1)}_t, \sigma^{(2)}_t$, which for
instance equals $1.3\times \int_0^T (1+\rho_t^2)(\sigma_t^{(1)})^4
dt$ if $(\sigma_t^{(1)})^2 \equiv (\sigma_t^{(2)})^2$, while here we
show that a CLT holds for the threshold estimator even in presence
of (finite activity) jumps, the asymptotic standard error being
$\int_0^T (1+\rho_t^2)(\sigma_t^{(1)})^2(\sigma_t^{(2)})^2 dt$, for
any $\rho_t, \sigma^{(1)}_t, \sigma^{(2)}_t$, and which is less than
the error of the BPC at least when
 $(\sigma_t^{(1)})^2 \equiv (\sigma_t^{(2)})^2$.
The bipower covariation test of \cite{BarShe04BPC} has been
discussed by \cite{BolLawTau07}, where the Authors show that, when
dealing with large portfolios, it is necessary to use a different
global cross-variation index to
get reliable results.

A CLT using multipowers for a bivariate process and in presence of
jumps is given by \cite{JacTor07}. More precisely, regarding the
co-jumps, they consider the quantity $ \hat{B}_n:=\sum_{i=1}^n
(\Delta_i X^{(1)})^2(\Delta_i X^{(2)})^2$, which is an estimate of
$B:=\sum_{s\leq T} (\Delta X^{(1)}_s)^2(\Delta X^{(2)}_s)^2$, not
directly comparable with an estimate of the sum of the co-jumps
$\sum_{s\leq T} \Delta X^{(1)}_s\Delta X^{(2)}_s$ we give here.
Their goal is to give a test for the presence of co-jumps, so they
concentrate on $$ \phi^{(j)}_n=\frac{\sum_{i} (\tilde\Delta_i
X^{(1)})^2(\tilde\Delta_i X^{(2)})^2}{\sum_{i} (\Delta_i
X^{(1)})^2(\Delta_i X^{(2)})^2},$$ the quotient of two cross-power
variations of the bivariate $X$, computed for different lags $kh$
and $h$: $\tilde\Delta_i X^{(q)}:= X^{(q)}_{ikh}-
X^{(q)}_{(i-1)kh}$, $\Delta_i X^{(q)}:= X^{(q)}_{ih}-
X^{(q)}_{(i-1)h}$. They reach that $ \phi^{(j)}_n\ri 1$ as $n\ri
\infty$, on the space $\Omega^{(j)}$ where some co-jumps occur,
and they prove a CLT for $ \phi^{(j)}_n$ in restriction to
$\Omega^{(j)}$.  We remark that to compute an estimator of the
conditional asymptotic variance of  $ \phi^{(j)}_n$ they in fact
use the threshold
technique when the volatilities are stochastic and are allowed to co-jump with the respective state variables.\\


An outline of the paper is as follows. In section 2 we illustrate
the framework; in section 3 we deal with the case where each
component $J^{(q)}$ of $X^{(q)}$ has finite activity of jump. We
show that $\hat{IC}_{T,n}$ is asymptotically Gaussian, so that it
is also consistent. We find a joint CLT allowing to estimate the
$\beta$s and the correlation coefficient of the continuous parts
of the two processes $X^{(q)}$, and we deal even with the case
where we dispose of non-synchronous observations. In section 4 we
deal with the more complex case where each $J^{(q)}$ can have an
infinite activity semimartingale jump component $\tilde
J_2^{(q)}$. We show that our estimator is still consistent. Since
the given theory asserts that we can asymptotically identify the
quantities of our interest, in section 5 we check on simulations
that in fact the finite sample performance of $\hat{IC}_{T,n}$ is
good even for time step between the
 observations large enough (five minutes) to avoid considering microstructure effects on the data, at least
for commonly used financial models with realistic choices of the
parameters. Section 6 concludes and section 7 contains all the
proofs and  technical details.

\section{Framework and notation}

Consider a filtered probability space
$(\Omega,\mathcal{F},(\mathcal{F}_{t})_{t\in [0,T]},P)$ where
$X^{(1)}=(X_{t}^{(1)})_{t \in [0,T]}$ and $X^{(2)}=(X_{t}^{(2)})_{ t
\in [0,T]}$ are two real Itô semimartingales  defined by
\begin{equation}\label{modello}
  X_{t}^{(q)} = \int_{0}^{t}a_{s}^{(q)}ds+\int_{0}^{t}\sigma_{s}^{(q)}
  dW_{s}^{(q)}+J_{t}^{(q)},\quad t\in [0,T],\quad q=1,2
\end{equation}
where

    \hspace{1cm}
      \begin{center}
     \begin{minipage}[c]{0.9\linewidth}
     \textbf{A1}. $W^{(1)}=(W^{(1)}_{t})_{t \in [0,T]}$ and  $W^{(2)}=(W^{(2)}_{t})_{t \in [0,T]}$
    are two correlated Wiener processes,
    with quadratic instantaneous covariation  given by\\
    $d<W^{(1)},W^{(2)}>_{t}$ $=\rho_{t} dt$, $t\in [0,T]$;
 \end{minipage}\end{center}
 we can write
    $$
    W^{(2)}_{t}=\rho_{t} W^{(1)}_{t}+\sqrt{1-\rho_{t}^{2}}\ W^{(3)}_{t},
    $$
    where $W^{(1)}$ and $W^{(3)}$ are independent standard Brownian motions.\\

   \begin{center}
      \begin{minipage}[c]{0.9\linewidth}\textbf{A2}.
     The diffusion stochastic coefficients $\sigma^{(q)}=(\sigma^{(q)}_{t})_{t\in
    [0,T]}$, $a^{(q)}=(a^{(q)}_{t})_{t\in [0,T]}$, $q=1,2$, and
    $\rho=(\rho_{t})_{t\in [0,T]}$ are  càdlàg  adapted processes.\\
    \end{minipage}\end{center}

\n As for the jump components $J^{(q)}$, in the next section we have FA jumps, i.e.
$$ J^{(q)}_{t}=
\sum_{k=1}^{N_{t}^{(q)}}\gamma_{\tau^{(q)}_{k}},\quad q=1,2,
$$
as specified with more detail below, where $N^{(q)}=(N_{t}^{(q)})_{t\in [0,T]}$ are counting processes
with $E[N^{(q)}_T]<\infty$.\\
More generally in section \ref{secIA} each $J^{(q)}$ is allowed to be any pure jump semimartingale
with possibly IA.\\

To begin with we assume to have equally spaced and synchronous
observations. The consistency results under not equally spaced but
synchronous observations are straightforward using Lemma
\ref{ModContD} and Theorem \ref{IndicatorNoJumps} below.
Generalization to not equally spaced and not synchronous
observations are dealt with later. Let, for each $n$,
$\pi_{n}=\{0=t_{0}^{(n)}<t_{1}^{(n)}<\cdot\cdot\cdot<t_{n}^{(n)}=T\}$
be a partition of $[0,T]$. For simplicity let us write $t_{j}$ in
place of $t_{j}^{(n)}$. Define $h:=t_{j}-t_{j-1}=\frac{T}{n},$ for
every $j=1,..,n$ and $n=1,2,..$. Note that $h\rightarrow 0$ if and
only if $n \rightarrow \infty$, so when computing the limits of
our interest we indifferently indicate one
of the two.\\

  \begin{center}\begin{minipage}[c]{0.9\linewidth}  \textbf{A3}. We choose a deterministic function, $h\mapsto
r(h)$, satisfying the following properties
$$
\lim_{h\rightarrow 0}
    r(h)=0, \quad \lim_{h\rightarrow
    0}\frac{h \log \frac{1}{h}}{r(h)}=0.
$$
\end{minipage}
\end{center}

\vm
\n We denote $r(h)$ by $r_{h}$, and, for each $q=1,2$,
$$ D^{(q)}_{t} = \int_{0}^{t}a^{(q)}_{s}ds+\int_{0}^{t}\sigma^{(q)}_{s}dW_{s}^{(q)},  $$
the Brownian semimartingale part of $X^{(q)}$.\\

\n As a consequence of the Paul L\'evy result 
about the modulus of continuity of the Brownian motion paths, we
can control how quickly the increments of the diffusion part of
each $\Delta_{j}X^{(q)}$ tend to zero. This is the key point to
understand when  $\Delta_{j}X^{(q)}$ is likely to contain some
jumps. More precisely, the Paul L\'evy law (see e.g.
\cite{KarShr99}, p.114, Theorem 9.25)
implies that
  \beq\label{mcBM}
  \mbox{a.s.} \quad  \lim\limits_{h\ri 0}\sup_{j\in \{ 1,..,n \} }\frac{|\DW|}{\sqrt{2h \log\frac{1}{h}} }
  \leq 1. \eeq
However the stochastic integral $\sigma.W$ is  a time changed
Brownian motion (\cite{RevYor01}, Theorems 1.9 and 1.10),
i.e. $ \Delta_{j}\left(\sigma.W\right)=
B_{IV_{t_{j}}}-B_{IV_{t_{j-1}}},$ where $B$ is a Brownian motion
and $IV_t$ is the integrated variance $\int_0^t\sigma^2_s ds$ up
to time $t$. Note that the increments of the drift part of $X$
tend to zero more quickly than $\sqrt{2h \log\frac{1}{h}}$ as
$h\ri 0$, so for $D^{(q)}$ we can reach a result similar to
(\ref{mcBM}), as soon as the boundedness of the paths of $a$ and
$\sigma$ is guaranteed (which is the case when they are càdlàg).
In fact
$$
\sup_{j=1..n} \frac{|\intIi a_sds+\intIi\sigma_s dW_s|}{\sqrt{2h
\log\frac{1}{h} }}\leq \sup_j\frac{|\intIi a_sds|}{\sqrt{2h
\log\frac{1}{h} }}+ \sup_j\frac{|\intIi\sigma_sdW_s|}{\sqrt{2h
\log\frac{1}{h} }}\leq$$
$$C(\omega)\sqrt{\frac{h}{\log\frac 1 h}} +
\sup_{j}\frac{|B_{IV_{t_{j}}}-B_{IV_{t_{j-1}}}|}{\sqrt{2\Delta_{j}
IV \log\frac{1}{\Delta_{j} IV }} }\
\sup_{j}\frac{\sqrt{2\Delta_{j}IV \log\frac{1}{\Delta_{j} IV
}}}{\sqrt{2M(\omega)h \log\frac{1}{M(\omega)h}} }
\sup_{j}\frac{\sqrt{2M(\omega) \log\frac{1}{M(\omega)h}} }{\sqrt{2
\log\frac{1}{h}} },
$$
where $C(\omega):=\sup_{s\in[0,T]} |a_s(\omega)|,
M(\omega):=\sup_{s\in[0,T]} |\sigma^2_s(\omega)|$. By
\cite{KarShr99} (Theorem 9.25) and the monotonicity of the
function $x\ln \frac 1 x$ it follows that
as $h\ri 0$, 
the right hand side has a limsup which is bounded by  $
\sqrt{M(\omega)}$, thus for sufficiently small $h$, even in the case of not equally spaced observations, the following
holds.


\begin{lemma}\label{ModContD} (\cite{mancini07}) Under \textbf{A2} we have that, given an arbitrary  partition
$\{t_0=0, t_1, ..., t_n=T\}$ of $[0,T],$ then for sufficiently
small $h:=\sup_{j=1..n} |t_{j}-t_{j-1}|$ we have  a.s.
$$
\sup_{j=1..n}\frac{|\Delta_{j}D^{(q)}|}{\sqrt{2h \log \frac{1}{h}}}
\leq K_{q}(\omega),\quad q=1,2,
$$
where $K_q(\omega) := \sqrt{M(\omega)} + 1$ are finite random variables.
\end{lemma}

Last result implies that if $\DXq>\rh$ and $\rh$ is, for small
$h$, larger than $2 K^2_q h\log\frac{1}{h}$ (as it is, under {\bf
A3}), then we have $\DXq>2 K^{2}_q h\log\frac{1}{h}$, and it is
not likely that $\DX$ coincides with the increment of a Brownian
semimartingale, while it is likely that some jumps occurred within
$\Ii$ and made $|\DX|$ large.

Application of Lemma \ref{ModContD} gives us the main tool for the
construction of our estimators in the next section.

{\bf Notation.}

$\bullet$ For any semimartingale $Z$,  $  \Delta Z_{s} =
Z_{s}-Z_{s-}$ denotes the size of the jump of $Z$ at time $s$, while
$ \Delta_j Z = Z_{t_j}-Z_{t_{j-1}}$ denotes the increment of process
$Z$ in the time interval $\Ii$


$\bullet$ $IC_t=\int_0^t\rho_{s}\sigma^{(1)}_s\sigma^{(2)}_s ds$
denotes the integrated covariation up to time $t$,\\
$\hat{IC}_{t, n}= \sum_{j=1..n: \ t_j\leq t}
\Delta_{j}X^{(1)}1_{\{(\Delta_{j}X^{(1)})^{2}\leq
r_h\}}\Delta_{j}X^{(2)}1_{\{(\Delta_{j}X^{(2)})^{2}\leq r_h\}}, \
h=T/n,$ is its threshold estimator

$\bullet$  $IV^{(q)}_t=\int_0^t(\sigma^{(q)}_s)^2 ds$ denotes the
integrated variance of process $X^{(q)}$, q=1,2, up to time $t$
and $\hat{IV}^{(q)}_{t, n}= \sum_{j=1..n: \ t_j\leq t}
(\Delta_{j}X^{(q)})^21_{\{(\Delta_{j}X^{(q)})^{2}\leq r_h\}}$ is
its threshold estimator

$\bullet$ sometimes
$\Delta_{j}X^{(q)}1_{\{(\Delta_{j}X^{(q)})^{2}\leq r_h\}}$ is
indicated briefly with $\Delta_{j}X^{(q)}_\star$

$\bullet$ sometimes we write Plim to indicate the limit in
probability. $\stackrel{st}{\ri}$ indicates stable convergence in
law of processes. See \cite{JacShi03}, ch. 8, sec. 5c, for the
definition and properties of stable convergence in law, and
\cite{JacDisp07}
 for further statement of useful properties


\section{Finite activity jumps: consistency and central limit theorem}\label{secFA}

In this section we assume that $J^{(q)}$ is any  FA jump process: 
for each $q=1,2$, 
$$
J^{(q)}_{t}=\int_{0}^{t}\gamma_{s}^{(q)}dN_{s}^{(q)}=\sum_{k=1}^{N_{t}^{(q)}}\gamma_{\tau^{(q)}_{k}},
$$
where $N^{(q)}=(N_{t}^{(q)})_{t\in [0,T]}$ is a counting process
with $E[N^{(q)}_T]<\infty$, $\{\tau^{(q)}_{k},\ k=1,...,N_{T}^{(q)}\}$ denote the
instants of jump of $J^{(q)}$ and $\gamma_{\tau^{(q)}_{k}}$ denote
the sizes $\Delta J^{(q)}_{t}$ of the jumps occurred at
$\tau^{(q)}_{k}$. Denote
$$\underline{\gamma}^{(q)}= \min_{k=1,...,N_{T}^{(q)}} |\gamma_{\tau^{(q)}_{k}}|.$$

  \begin{center}
     \begin{minipage}[c]{0.9\linewidth} \textbf{A4}.
\quad Assume $ E[N^{(q)}_T]<\infty$ and
$P(\gamma_{\tau^{(q)}_{k}}=0)=0,\quad \forall \
k=1,...,N_{T}^{(q)},\ q=1,2. $
 \end{minipage}
 \end{center}

\vspace{0.5cm}
\begin{remark}
Condition \textbf{A4} implies that a.s. $\underline{\gamma}^{(q)}>0$.
\end{remark}

\begin{example}
{\rm If  $J^{(q)}$ are FA
 L\'evy processes, then they are of
 compound Poisson type (\cite{ConTan04}, Proposition 3.3, section 3.2):
 $N^{(q)}$ are simple Poisson processes with constant intensities $\lambda^{(q)}$ and for each $q$ the random
 variables $\gamma_{\tau^{(q)}_{k}}$ are i.i.d., for $k=1,...,N_{T}^{(q)},$ are independent on $N^{(q)}$
 and satisfy condition \textbf{A4}.}\end{example}

We remark that the consistency and CLT we reach in this section
are valid in presence of  general finite activity jump processes,
in that we do not need any assumptions on the law of the jump
sizes, or of the counting processes $N^{(q)}$, nor any assumption
of independence. We do not even need that $J^{(q)}$ are FA jumping
semimartingales, we only need that {\bf A4} holds, which is true
if $J^{(q)}$ are (FA jumping) semimartingales.

 Now we construct our threshold estimators.

\begin{definition} We define for $r,l\in \N$
$$
v^{(n)}_{r,l}(X^{(1)},X^{(2)})_{t}
=h^{1-\frac{r+l}{2}}\sum_{j: t_j\leq t}(\Delta_{j}X^{(1)})^{r}(\Delta_{j}X^{(2)})^{l},
$$
$$
w^{(n)}(X^{(1)},X^{(2)})_{t}=h^{-1}\sum_{j: t_{j+1}\leq t}
\prod_{i=0}^{1}\Delta_{j+i}X^{(1)}\prod_{i=0}^{1}\Delta_{j+i}X^{(2)}.
$$
and their analogous \emph{threshold} versions
$$
\tilde{v}^{(n)}_{r,l}(X^{(1)},X^{(2)})_{t}
=h^{1-\frac{r+l}{2}}\sum_{j: t_j\leq t}(\Delta_{j}X^{(1)})^{r}1_{\{(\Delta_{j}X^{(1)})^{2}\leq
r_{h}\}}(\Delta_{j}X^{(2)})^{l}1_{\{(\Delta_{j}X^{(2)})^{2}\leq
r_{h}\}},
$$
$$
\tilde{w}^{(n)}(X^{(1)}\!\!,X^{(2)})_{t}\!=\!h^{-1}\!\sum_{j: t_{j+1}\leq t}
\prod_{i=0}^{1}\Delta_{j+i}X^{(1)}1_{\{(\Delta_{j+i}X^{(1)})^{2}\leq
r_{h}\}}\prod_{i=0}^{1}\Delta_{j+i}X^{(2)}1_{\{(\Delta_{j+i}X^{(2)})^{2}\leq
r_{h}\}}.
$$
\end{definition}

\n $v^{(n)}_{r,l}(X^{(1)},X^{(2)})_{T}$ and
$w^{(n)}(X^{(1)},X^{(2)})_{T}$ are used in \cite{BarShe04}
 to estimate $IC_T$ in the case where $X^{(q)}$ are diffusion processes. $\tilde
v^{(n)}_{r,l}(X^{(1)},X^{(2)})_{T}$ and $\tilde
w^{(n)}(X^{(1)},X^{(2)})_{T}$ are modified versions
for the case of jump-diffusion processes: 
 by Theorem \ref{IndicatorNoJumps} they exclude from the sums the
terms containing  jumps. Note that $\tilde v^{(n)}_{1,1}(X^{(1)},X^{(2)})_{t}= \hat{IC}_{t,n},$ for all $t\in[0,T]$.

\n In view of the practical application of our estimator we are now
interested in the speed of convergence of $\hat{IC}_{T,n}$. We in fact reach even more. The first main
result of this section is  a joint central limit theorem
 for the threshold estimators 
$$\left(\begin{array}{ll}
\hat{IV}^{(1)} & \hat{IC}\\
\hat{IC} & \hat{IV}^{(2)}
\end{array}\right)$$
which implies that in presence of finite activity jumps
$\hat{IC}_{T,n}$ converges to $IC$ at speed $\sqrt h$, $h=T/n$,
and it allows to give estimators of standard dependence measures
between the diffusion parts $D^{(q)}$ of our processes $X^{(q)}$,
such as the {\it realized diffusion regression} coefficients up to
time $t$
$$
\beta^{(1,2)}_t:= \frac{ IC_t }{IV^{(2)}_t}, \quad
\beta^{(2,1)}_t:= \frac{ IC_t }{IV^{(1)}_t}$$
and the
{\it realized diffusion correlation}
$$ \rho^{(1,2)}_t:= \frac{ IC_t }{ \sqrt{ IV^{(1)}_t IV^{(2)}_t  }  }.$$

\begin{theorem} [Joint CLT, FA jumps] \label{JointCv} Under assumptions from \textbf{A1} to \textbf{A4},
with $h=T/n,$ we have, as $h\ri 0$,
$$h^{-1/2}\left(\begin{array}{cc}
\hat{IV}^{(1)}_n- {IV}^{(1)} & \hat{IC}_n- IC\\
\hat{IC}_n -IC & \hat{IV}^{(2)}_n- {IV}^{(2)}
\end{array}\right) \stackrel{st}{\ri}
\frac{1}{\sqrt 2}\left(\begin{array}{cc}
2Z_{11} & Z_{12} + Z_{21}\\
Z_{12} + Z_{21} & 2Z_{22}
\end{array}\right),$$
where ${\bf Z}$ is the $2\times 2$ process with components
$$Z_{11, t}:=  \int_0^t (\sigma^{(1)}_s)^2dB_{11 s} $$
\beq\label{AsLawnelJointCLT} \begin{array}{c}
 Z_{12, t}:= \int_0^t \rho_s\sigma^{(1)}_s\sigma^{(2)}_sdB_{11 s} +
  \int_0^t \sqrt{1-\rho_s^2}\ \sigma^{(1)}_s\sigma^{(2)}_sdB_{12s }\\
Z_{21, t}:= \int_0^t \rho_s\sigma^{(1)}_s\sigma^{(2)}_sdB_{11 s} +
\int_0^t \sqrt{1-\rho_s^2}\ \sigma^{(1)}_s\sigma^{(2)}_sdB_{21 s} \end{array}\eeq
$$ \!\!\!\!
Z_{22, t}:=\int_0^t \rho^2_s(\sigma^{(2)}_s)^2dB_{11 s} + \int_0^t \rho_s\sqrt{1-\rho_s^2}\ (\sigma^{(2)}_s)^2
                         \Big(dB_{12 s}+dB_{21 s} \Big)+ \int_0^t (1-\rho_s^2)(\sigma^{(2)}_s)^2 dB_{21 s}$$

\n and ${\bf B}$ is a $2\times 2$-dimensional  standard Brownian motion independent on the filtered probability
space
$(\Omega, {\cal F}, ({\cal F}_t)_{t\in[0,T]}, P)$ where our model
is defined.

\end{theorem}

 Note that the result  for $\hat{IV}^{(q)}_n$ is consistent with \cite{mancini07},
 since $$Var(\sqrt 2 Z_{qq, T})= 2\int_0^T \rho^4_s(\sigma^{(q)}_s)^4ds +
                     4\int_0^T \rho^2_s (1-\rho_s^2)(\sigma^{(q)}_s)^4 ds
                     + 2\int_0^T (1-\rho_s^2)^2(\sigma^{(q)}_s)^4 ds$$
                     $$= 2\int_0^T (\sigma^{(q)}_s)^4ds.$$

\begin{corollary} [Consistency, FA jumps]\label{hatICconsistFA} Under {\bf A1} to {\bf A4}, as $n\ri \infty$, for all
$t\in[0,T]$
$$\hat{IC}_{t, n}\stackrel{P}{\ri}  IC_t,$$
if {\it a.s} $IV_t^{j}\neq 0$ then
$$\hat{\beta}^{(i,j)}_{t, n} := \frac{ \hat{IC}_{t, n} }{\hat{IV}^{(j)}_{t,n}}\stackrel{P}{\ri} \beta^{(i,j)}_t, \quad (i,j)= (1,2), (2,1)$$
if {\it a.s} $IV_t^{1}IV_t^{2}\neq 0$ then
$$ \hat{\rho}^{(1,2)}_{t, n}:= \frac{ \hat{IC}_{t,n} }{ \sqrt{ \hat{IV}^{(1)}_{t,n} \hat{IV}^{(2)}_{t, n}   }  }
\stackrel{P}{\ri} \rho^{(1,2)}_t.$$
\end{corollary}


\begin{corollary} [Speed of convergence of $\beta$s and $\rho$, FA jumps] \label{speedCvBetasRho}
If a.s. ${IV}^{(j)}_t\not = 0$ for all $t\in [0,T]$ we have, for $(i,j)= (1,2)$ or $(2,1)$,
$$h^{-1/2}\left(\hat{\beta}^{(i,j)}_n - \beta^{(i,j)}\right)
 \stackrel{st}{\ri} \frac{Z_{12}+ Z_{21}}{\sqrt 2 IV^{(j)}} + \sqrt 2 Z_{jj}\frac{IC}{(IV^{(j)})^2}.$$
If a.s. ${IV}^{(1)}_t{IV}^{(2)}_t\not= 0$ for all $t\in [0,T]$ we have
$$h^{-1/2}\left(\hat{\rho}^{(1,2)}_n - \frac{IC}{\sqrt{IV^{(1)}IV^{(2)}}}\right)
  \stackrel{st}{\ri}$$
  $$\frac{Z_{12}+ Z_{21}}{\sqrt {2 IV^{(1)}IV^{(2)}}}- \frac{ Z_{22} IC}{\sqrt{2\  IV^{(1)}} (IV^{(2)})^{3/2} }-
\frac{Z_{11} IC }{\sqrt{2\  IV^{(2)}} (IV^{(1)})^{3/2} }.$$
\end{corollary}

\vspace{0.5cm} The following proposition allows us to give a CLT
for the standardized version of the estimation error
$\hat{IC}_{T,n}- IC_T$. Note that the asymptotic variance of
$h^{-1/2}(\hat{IC}_{T,n} -IC_T)$, by Theorem \ref{JointCv}, is
given by $(Var(Z_{12,T}+Z_{21,T}))/2=
\int_{0}^{T}(1+\rho_{t}^{2})(\sigma_{t}^{(1)})^{2}(\sigma_{t}^{(2)})^{2}dt.$

\begin{proposition} [Estimate of the standard error for $\hat{IC}_{t,n}$, FA jumps] \label{DenominPerCLTFA}
Under assumptions \textbf{A1} to \textbf{A4}
we have, for all $t\in[0,T]$,
$$
\tilde{v}^{(n)}_{2,2}(X^{(1)},X^{(2)})_{t}-\tilde{w}^{(n)}(X^{(1)},X^{(2)})_{t}\stackrel{P}{\longrightarrow}
\int_{0}^{t}(1+\rho_{s}^{2})(\sigma_{s}^{(1)})^{2}(\sigma_{s}^{(2)})^{2}ds.
$$
\end{proposition}

\vspace{0.5cm} \n We now are ready to present  the central limit theorem for the
standardized estimation error.

\begin{corollary} [CLT for the standardized version of  $\hat{IC}_{t,n}-IC_t$, FA jumps] \label{CLTFiniteActivity}
Under \textbf{A1} to \textbf{A4}, if
a.s.  $\int_{0}^{t}(1+\rho_{s}^{2})(\sigma_{s}^{(1)})^{2}(\sigma_{s}^{(2)})^{2}ds\not=0$
 we have
$$
\frac{\hat{IC}_{t,n}-IC_t}
{\sqrt h\sqrt{\tilde{v}^{(n)}_{2,2}(X^{(1)},X^{(2)})_{t}-\tilde{w}^{(n)}(X^{(1)},X^{(2)})_{t}}}
\stackrel{d}{\longrightarrow} {\cal N},
$$
where ${\cal N}$ denotes a standard Gaussian random variable.
\end{corollary}



\begin{remark} [Estimate of the co-jumps] {\rm By Corollary \ref{hatICconsistFA}, clearly we have an estimate of the
sum of the co-jumps up to $T$ simply subtracting $\hat{IC}_{T,n}$  from the quadratic covariation estimator:
$$ \sum_{j=1}^{n}\Delta_{j}X^{(1)}\Delta_{j}X^{(2)}-
\hat{IC}_{T,n}
\stackrel{P}\longrightarrow \sum_{0\leq s\leq T} \Delta
J_s^{(1)}\Delta J_s^{(2)},
$$
as $n\ri \infty$. Analogously we can obtain an estimator of the sum of the co-jumps up to each time $t\in[0,T]$.\\
An estimate of each $ \Delta J_s^{(1)}\Delta J_s^{(2)}$, with $s\in[0,T]$, is obtained using
\beq\label{StimSingleCojleq}\Delta_{j}X^{(1)}\Delta_{j}X^{(2)}-\Delta_{j}X^{(1)}1_{\{(\Delta_{j}X^{(1)})^{2}\leq r_{h}\}}
\Delta_{j}X^{(2)}1_{\{(\Delta_{j}X^{(2)})^{2}\leq r_{h}\}}, \eeq with
$j$ such that $s\in ]t_{j-1}, t_j].$ Alternatively, as we
consider one single term, and not the sum of $n$ terms, even
\beq\label{StimSingleCojmag}
\Delta_{j}X^{(1)}1_{\{(\Delta_{j}X^{(1)})^{2}> r_{h}\}}
\Delta_{j}X^{(2)}1_{\{(\Delta_{j}X^{(2)})^{2}> r_{h}\}}\eeq or
\beq\label{StimSingleCojnoI}\Delta_{j}X^{(1)}\Delta_{j}X^{(2)}\eeq estimate the co-jump $ \Delta
J_s^{(1)}\Delta J_s^{(2)}$, with $s\in ]t_{j-1}, t_j],$ since
%
%
$|\Delta_{j}X^{(q)}1_{\{(\Delta_{j}X^{(q)})^{2}\leq r_{h}\}}
\Delta_{j}X^{(\ell)}|$ $\leq 2\sqrh \sup_{s\in[0,T]} |X^{(\ell)}|$, $q=1,2$, $\ell=3-q$, and
$|\Delta_{j}X^{(1)}|1_{\{(\Delta_{j}X^{(1)})^{2}\leq r_{h}\}}
|\Delta_{j}X^{(2)}|\times $ $1_{\{(\Delta_{j}X^{(2)})^{2}\leq r_{h}\}}\leq \rh$
 tend to
zero in probability as $h\ri 0,$ by the pathwise boundedness of each $X^{(\ell)}$ on $[0,T]$.
However, as we show in section 5, estimator (\ref{StimSingleCojmag}) has the
best finite sample properties in the simulations of Model 1 having FA jumps. \qed
 }\end{remark}


\begin{remark} {\bf Finite sample performance and microstructure noises}.
{\rm Our theoretic results allow to estimate $IC_T$ and the co-jumps asymptotically for $h\ri 0$, while
in practice for very small values of $h$ financial time series are affected by microstructure noises which
introduce a bias which is larger as $h$ is smaller.
In section 5 we implement our estimators of the integrated covariance and of the
co-jumps on simulations of realistic financial time series and we find that they
have good performance already with temporal mesh $h$
corresponding to five minutes, a time lag at which prices are not usually affected by microstructure noises
 (\cite{BS07Variation}).
However we remark that when iid microstructure noises
contaminate the observations of each asset price $X^{(q)}$, the threshold estimator
rules out even the noises, similarly as
it rules out the contribution of the jumps (\cite{Mancini08}).}
\end{remark}

\vspace{0.5cm}
\begin{remark} {\bf Asynchronous observations.} {\rm
It is known that the problem of the estimation of the covariation
among two assets undergoes the so called Epps effect, i.e. in the
empirical applications the estimator tends to zero as the step of
observation $h$ tends to zero. The asynchronicity among the
observations of $X^{(1)}$ and $X^{(2)}$ is considered one of the
possible causes (\cite{Ren03Epps}; \cite{BS07Variation}, section 2.10.2).  In fact some
Authors have tackled the problem of reaching a consistent
estimator of the covariation even when data are asynchronous
and $h\ri 0$, under the assumption of 
Brownian semimartingale models (in \cite{HayYos05} the estimator is introduced,
however we refer to
\cite{HayKus08} where the observation times are allowed to be dependent on $X^{(q)}$).\\
At the time scale of five minutes
the Epps effect probably does not affect our estimate of $IC_T$. However even in presence of
this microstructure-type noise (for smaller $h$) it is possible to make our estimator
correctly converge to the integrated covariation, as detailed
below.}
\end{remark}

Assume we dispose of two records  $\{D^{(1)}_{\tau^{(n)}_{0}},
D^{(1)}_{\tau^{(n)}_{1}}, ... D^{(1)}_{\tau^{(n)}_{m^{(n)}}}\}$, $\{D^{(2)}_{\nu^{(n)}_{0}},
D^{(2)}_{\nu^{(n)}_{1}}, ... D^{(2)}_{\nu^{(n)}_{k^{(n)}}}\}$, of observations of
two Brownian semimartingales $D^{(1)}$ and $D^{(2)}$, with the two stochastic partitions
 $0=\tau^{(n)}_{0}< \tau^{(n)}_{1}< ... \tau^{(n)}_{m^{(n)}}$ and  $0=\nu^{(n)}_{0}< \nu^{(n)}_{1}< ...
\nu^{(n)}_{k^{(n)}}$ spanned on $[0,T]$. For simplicity let us write
$\nu_{i}$ and $\tau_{j}$ in place of $\nu^{(n)}_{i}$ and $
\tau^{(n)}_{j}$. The idea of Hayashi and Kusuoka is to select only
some of the cross variations
$(D^{(1)}_{\tau_{j}}-D^{(1)}_{\tau_{j-1}})
(D^{(2)}_{\nu_{i}}-D^{(2)}_{\nu_{i-1}})$, in order to estimate the
covariation, and precisely  the ones for which there is an
intersection between the time intervals $]\tau_{j-1}, \tau_{j}]$ and
$]\nu_{i-1}, \nu_{i} ]$.\\
We show here that using their result (\cite{HayKus08}, Corollary 2.2) we in fact
reach the same kind of consistency in the case of asynchronous observations even
in presence of finite activity jumps. The idea is very
simple: we first eliminate the jumps, using threshold technique, and
then apply the Hayashi and Kusuoka  estimator to the estimated
continuous components $\hat D^{(q)}$. Recall that
$$ (X^{(q)}_{a}-X^{(q)}_{b})_\star = (X^{(q)}_{a}-X^{(q)}_{b}) 1_{\{(X^{(q)}_{a}-X^{(q)}_{b})^2 \leq \rh\}},$$
for any two time instants $a$ and $b$, $h:=\sup_{j=1..m^{(n)}} (\tau_{j}- \tau_{j-1}) \vee \sup_{i= 1..k^{(n)}}
(\nu_{i}- \nu_{i-1}).$ We in fact have the following

\begin{theorem} [Asynchronous observations]\label{asyncObs}
Let  {\bf A1} to {\bf A4}
hold, $0=\tau_{0}< \tau_{1}< ... < \tau_{m^{(n)}}$,
 $0=\nu_{0}< \nu_{1}< ... < \nu_{k^{(n)}}$ be two sequences of  stopping times such that
$\tau_{m^{(n)}}\uparrow T$, $\nu_{k^{(n)}}\uparrow T$
a.s., as $n\ri \infty$ then
$$ \sum_{j=1..m^{(n)}, i= 1..k^{(n)}}
(X^{(1)}_{\tau_{j}}-X^{(1)}_{\tau_{j-1}})_\star\
(X^{(2)}_{\nu_{i}}-X^{(2)}_{\nu_{i-1}})_\star\
1_{\{]\tau_{j-1}, \tau_{j}]\cap]\nu_{i-1}, \nu_{i} ] \neq
\emptyset \}} \stackrel{P}{\ri} IC_T,
$$
as $h:=\sup_{j=1..m^{(n)}} (\tau_{j}- \tau_{j-1}) \vee \sup_{i= 1..k^{(n)}}
(\nu_{i}- \nu_{i-1})\stackrel{P}{\ri} 0. $
\end{theorem}

\section{Infinite activity jumps: consistency}\label{secIA}

In this section we allow the jump components of processes $X^{(q)}$ to have infinite activity, so we are here in the
case where $X^{(q)}$ are general Itô semimartingales.
Any unidimensional Itô semimartingale has a representation as in (\ref{modello}) with each $J^{(q)}$ decomposed as  \begin{equation}\label{decompJ}
\begin{array}{c}
J^{(q)}=J_{1}^{(q)}+\tilde{J}_{2}^{(q)},\\
 J^{(q)}_{1 t}(\omega)=\int_{0}^{t}\int_{|\gamma^{(q)}(\omega,t,x)|>1} \gamma^{(q)}(\omega,t,x)
\underline{\mu}^{(q)}(\omega, dx,dt),\\
$$\tilde{J}_{2t}^{(q)}(\omega)=\int_{0}^{t}\int_{ |\gamma^{(q)}(\omega,t,x)|\leq 1}\gamma^{(q)}(\omega,t,x)
\tilde{\underline{\mu}}^{(q)}(\omega,dx, ds),
\end{array}
\end{equation}
where $\underline{\mu}^{(q)}$ is the Poisson random measure of the
jumps of $J^{(q)}$, $\tilde{\underline{\mu}}^{(q)}(\omega, dx,
ds)=\mu^{(q)}(\omega, dx, ds)- \underline{\nu}^{(q)}(\omega, dx,
ds)$ is its compensated measure, $\underline{\nu}^{(q)}(\omega, dx,
ds)= dx\times ds$, the coefficients $a^{(q)}, \sigma^{(q)},
\gamma^{(q)}$ are predictable and
 $\int 1\wedge (\gamma^{(q)})^2(\omega, t,x) dx$ is a.s. finite (see \cite{JacDisp07}, pp.3,4;
 \cite{Jac08}, (2.11)).\\
Conditions {\bf A2} and {\bf A4'} below guarantee local boundedness properties of  such coefficients.

    \begin{center}
     \begin{minipage}[c]{0.9\linewidth} \textbf{A4'}.
     $\int 1\wedge (\gamma^{(q)})^2(\omega, t,x) dx$ is  locally bounded.
 \end{minipage}
 \end{center}

\n For each $q=1,2$, $J_1^{(q)}$ is a finite activity jump process
of type $ J^{(q)}_{1
t}=\sum_{k=1}^{N_{t}^{(q)}}\gamma_{\tau^{(q)}_{k}}, $ as in
section \ref{secFA}, where $E[N^{(q)}_T]<\infty$ is equivalent to
$E[\int_{|\gamma^{(q)}|>1} dx]<\infty$ and now the sizes
$|\gamma_{\tau^{(q)}_{k}}|$ are all larger than 1;
on the contrary $\tilde{J}_{2}^{(q)}$ accounts for
the infinite activity jumps of $J^{(q)}$, since generally $\int_{|\gamma^{(q)}|\leq
1} dx =+\infty$.
$\tilde{J}_{2}^{(q)}$ is a
compensated sum of  jumps, where each jump is  bounded in absolute
value by 1. Therefore, for each $q=1,2$, $J^{(q)}_{1}$  accounts for the
"large" 
and rare jumps  of $X^{(q)}$, while $\tilde J_2^{(q)}$ accounts
for the  frequent and small jumps.

\begin{example}\label{RemCasoLevyJumps}
{\rm If one of the two processes $J^{(q)}$ is a pure jump L\'evy process, it is always possible to decompose it as in
(\ref{decompJ}) with $\gamma^{(q)}(\omega, t, x)\equiv x$ but
$\underline{\nu}^{(q)}(\omega, dx, ds)\equiv \nu^{(q)}(dx)\times ds$, 
where $\nu^{(q)}$
is the Lévy measure of $J^{(q)}$ and is a deterministic $\sigma$-finite measure such that
$\int_{\R} 1\wedge x^2 \ \nu^{(q)}(dx)<\infty$ but generally such that $\int_{|x|\leq 1} \nu^{(q)}(dx)=+\infty$.
}\end{example}

 We prove that $\hat{IC}_{t,n}$ is still a consistent estimator of
$IC_t$, for all $t\in[0,T]$. For ease of notation we only consider
$IC$ up to time $T$ and evenly spaced synchronous observations.
Not evenly spaced but synchronous observations (with $h=
\sup_j|t_j-t_{j-1}|$) and arbitrary $t \in [0,T]$ are
straightforward. As a consequence the same estimators of the
co-jumps, presented in the previous section, are  consistent
even in the present framework.\\
As for the speed of convergence of
$\hat{IC}_{T, n}$, in the presence of
 infinite activity jump components, things are more complicated in that such a speed
 is determined both by the dependence
structure between $\tilde J_2^{(1)}$, $\tilde J_2^{(2)}$ and by the amount of jump activity of each
  $\tilde J_2^{(q)}$.
In \cite{GobMan07} we consider two Lévy infinite activity jump
components $\tilde J^{(1)}_2$ and $\tilde J^{(2)}_2$ with a
dependence structure described by a Lévy copula. We find that,
when $\tilde J^{(1)}_2$ and $\tilde J^{(2)}_2$ do depend, the
speed is still $\sqrt h$ only when the activity of jump of at
least one process is moderate (Blumenthal-Getoor index smaller
than 1), otherwise the speed is less than $\sqrt h$.

We now state the main result in presence of infinite activity jumps.

\begin{theorem}[Consistency in presence of IA jumps, synchronous observations]\label{consistencyIA}
Let $(X^{(1)}_{t})_{t\in [0,T]}$ and
$(X^{(2)}_{t})_{t\in [0,T]}$ be two processes of the form (\ref{modello}). Assume
\textbf{A1}, \textbf{A2}, \textbf{A3} and \textbf{A4'}. Then
$$
\hat{IC}_{T,n}\stackrel{P}{\longrightarrow} IC_T,
$$
as $n\rightarrow \infty$.
\end{theorem}


\vspace{-0.5cm}
\begin{remark} [Estimate of the co-jumps]
{\rm
Even in this framework of infinite activity jumps, as a consequence
of Theorem \ref{consistencyIA}, the sum of the co-jumps up to $T$ is estimated
by
$$\sum_{j=1}^n \Delta_{j}X^{(1)}\Delta_{j}X^{(2)}-\sum_{j=1}^n \Delta_{j}X^{(1)}1_{\{(\Delta_{j}X^{(1)})^{2}\leq r_{h}\}}
\Delta_{j}X^{(2)}1_{\{(\Delta_{j}X^{(2)})^{2}\leq r_{h}\}}.$$
Estimates of a single co-jump $\Delta J^{(1)}_s\Delta J^{(2)}_s$, with $s\in]t_{j-1}, t_j]$,
exactly as in section \ref{secFA}, are  given by (\ref{StimSingleCojleq}), (\ref{StimSingleCojmag})
 or (\ref{StimSingleCojnoI}). Simulations in section \ref{simulsStimaSinglecoj} show that
 for Model 2, with IA jumps, in fact estimator (\ref{StimSingleCojmag}) is a little bit more biased than (\ref{StimSingleCojnoI})
 but is still acceptable. Note that since in sec. 5 each $J^{(q)}$ has infinite activity and
 $J^{(2)}=\rho_J J^{(1)}+\sqrt{1-\rho_J^{2}}J^{(3)}$, each
$\Ii$  contains an infinite number of co-jump instants. \qed
 }\end{remark}

\section{Implementation}

\subsection{Choice of the threshold}\label{OptThr}
Our estimators depend on the threshold function $\rh$. In this
section we check on simulations how the results are sensitive to the
choice of $\rh$ in a given class. This is only an informal and
necessarily limited investigation. Formal study of methods for optimal
threshold selection in a given model is object of further
research.\\
In principle there are many functions $\rh$ satisfying conditions
{\bf A3}. However on simulations we find that the choice of $\rh$
within the family of powers of $h$, $\rh=ch^\beta$, with $c$ a
constant and $\beta$ a power in $]0,1[$, seems to be sufficiently good.\\
We simulate two kind of models: Model 1, proposed in
\cite{HuaTau05}, where each $X^{(q)}$ has stochastic volatility and
a  FA Compound Poisson jump part and Model 2, proposed in
\cite{CarGemMadYor02}, where each $X^{(q)}$ has constant volatility
and IA jumps,
as described in Table A. For Model 1 the parameters of the
univariate $X^{(q)}$ are taken from \cite{HuaTau05}. A path of each
$\sigma$ varies most between 0.013 and 0.019 in a day. For Model 2
the parameters of the univariate $X^{(q)}$ are taken from Table 2 of
\cite{CarGemMadYor02} for GE and HWP stocks. Note that the parameter
 $Y$ is not significantly different from zero for the two
 considered stocks, so that the CGMY process can be reduced to the VG
 process.

The VG process is characterized by three parameters $\kappa$,
$\theta$ and $\varsigma$. It is obtained by evaluating a Brownian
motion with drift, $\theta t + \varsigma B_t$, at a random time
$G_t$ given by a gamma process, a L\'evy process whose lag $h$
increments $G_{t+h}-G_t$ are distributed as Gamma r.v.s with mean
$h$ and variance $h\kappa$. It turns out that the VG process is pure
jump and has infinite, but moderate, activity (it is a process with
finite variation).

To effectively introduce non zero co-jumps, in each model the jump
component $J^{(2)}$ of $X^{(2)}$ is correlated with $J^{(1)}$ of
$X^{(1)}$ in the following way: we generate $J^{(1)}$ and an
independent $J^{(3)}$ with parameters as in Table A, then
$J^{(2)}=\rho_J J^{(1)}+\sqrt{1-\rho_J^{2}}J^{(3)}$. The simulation
of the model paths has been made using the  Euler scheme with
increments of 1 second, then we have taken the five minutes synchronous
returns and constructed our daily threshold estimator
$\hat{IC}_{T,n}$. We simulated 3000 bivariate paths.

For each model we implement the estimator of $IC$ as $\rh$ varies.
Figures 1-2 show how the mean relative bias in percentage form
$$100\frac{(\hat{IC}_{T,n}-IC)}{IC}$$
varies as $\beta$ varies in $]0,1[$ for $c=0.1,....,5.6$ with step
$0.5$ in Models 1 and 2, for $h$ fixed equal to five minutes ($n=84$
observations per day, time unit of measure $T$=1 day, $h=1/84$).
It is evident that the choice $c=0.1$ is the best one since in
presence of FA jumps ($\lambda^{(q)}=0.118$) it allows to decrease
the bias as $\beta$ increases. In fact in the case of IA jumps the
bias is much larger but $c=0.1$ allows to reach, for high $\beta$,
the lowest possible error. Figures 3-4 show the empirical
densities and the QQ-plots of the normalized bias
\begin{equation}\label{NB} \frac{\hat{IC}_{T,n}-IC_T} {\sqrt
h\sqrt{\tilde{v}^{(n)}_{2,2}(X^{(1)},X^{(2)})_{T}-\tilde{w}^{(n)}(X^{(1)},X^{(2)})_{T}}}
\end{equation}
when $r_h$ varies as before, for fixed $h$ equal to five minutes, for Model 1
with  $\lambda^{(q)}=0.014$. The same plots for Model 1 with
$\lambda^{(q)}=0.118$ and Model 2 are shown in Figures 5-6 and 7-8
respectively. We conclude that the best choice is
$r_h=0.1h^{0.99}$. As a further check in Figures 9-10 we made the
same plots for Model 1 with $\lambda^{(q)}=0$ and we found that
the choice of $r_h$ gives good results as well.

\subsection{Estimates of $IC$ and $\sum_{0\leq t\leq T}\Delta
J^{(1)}_{t}\Delta J^{(2)}_{t}$ on simulations} We report here the
performance of the estimators of $IC_T$ and of the sum
$\sum_{0\leq t\leq T}\Delta J^{(1)}_{t}\Delta J^{(2)}_{t}$ of the
co-jumps up to time $T$, where the threshold is the one selected
in the previous subsection. $T$ is kept fixed to one day, $h$ equals five
minutes.
 Figures 11-12-13 show the
histograms of $100\frac{(\hat{IC}_{T,n}-IC)}{IC}$ to check the
efficiency of $\hat{IC}_{T,n}$ for Model 1, $\lambda^{(q)}=0.014$,
Model 1, $\lambda^{(q)}=0.118$ and Model 2 respectively. The
relative summary statistics are shown in Tables 2 and 3.
$\hat{IC}_{T,n}$ has an acceptable performance in Model 1 and it is
biased in Model 2 but note that the estimation errors for $IC_T$ and
$\int_0^T (1+\rho_t^2)(\sigma_t^{(1)})^2(\sigma_t^{(2)})^2 dt$
compensate and give good empirical densities of the normalized bias
in Figure 7 for $c=0.1$ and $\beta=0.99$ and Table 1. Figures
14-15-16 show the histograms of the following relative bias in
percentage form
$$100\frac{(\sum_{j=1}^{n}\Delta_{j}X^{(1)}\Delta_{j}X^{(2)}-
\hat{IC}_{T,n})-\sum_{0\leq t\leq T}\Delta X^{(1)}_{t}\Delta
X^{(2)}_{t}}{\sum_{0\leq t\leq T}\Delta X^{(1)}_{t}\Delta
X^{(2)}_{t}}$$ for the sum of the co-jumps in Model 1,
$\lambda^{(q)}=0.014$, Model 1, $\lambda^{(q)}=0.118$ and Model 2.
Note that since in Model 2 each $J^{(q)}$ is a pure jump process with  IA in fact each movement
of $J^{(q)}$ is a jump, and each time that $J^{(1)}$ jumps even $J^{(2)}$ does by the way we correlated them, 
therefore the best we can do to reach the true $\sum_{0\leq t\leq
T}\Delta X^{(1)}_{t}\Delta X^{(2)}_{t}$ is to take the sum of the
cross-products of the one-second differences of processes $J^{(q)}$.
Tables 4 and 5 show the relative summary statistics. The performance
of the estimator of the sum of co-jumps is very good under Model 1
and a bit worse under Model 2. Under model 2 the estimate of
$\sum_{0\leq t\leq T}\Delta J^{(1)}_{t}\Delta J^{(2)}_{t}$ is much
better than the one of $IC_T$.




\subsection{Estimate of the single co-jumps}\label{simulsStimaSinglecoj}
Using the threshold function selected in subsection \ref{OptThr},
both for Model 1 and for Model 2
 we implement (\ref{StimSingleCojleq}), (\ref{StimSingleCojmag})
 and (\ref{StimSingleCojnoI})
to estimate each single co-jump, in order to check which
 is the most informative estimator for the single co-jumps. We consider 1 day time horizon and $h$ equal
 to five minutes.
 Figures 17-18 show the histograms of the
3000 values of $100\frac{\hat{JJ}_t - \Delta J^{(1)}_{t}\Delta
J^{(2)}_{t}}{\Delta J^{(1)}_{t}\Delta J^{(2)}_{t}}$ for each
estimator for both Model 1 ($\lambda_1=\lambda_2=0.118$) and Model
2, where we define by
 $\hat{JJ}_t$ (joint jumps) the estimate of $ \Delta J^{(1)}_{t}\Delta
 J^{(2)}_{t}$. Tables 6-7 report the relative summary statistics.
 We conclude that the most informative estimate of $\Delta
J^{(1)}_{t}\Delta J^{(2)}_{t}$ is (\ref{StimSingleCojmag}) for
Model 1 and (\ref{StimSingleCojnoI}) for Model 2, however, we find
that (\ref{StimSingleCojmag}) for Model 2 is still well
acceptable. We find that, anyway all three estimators show a good
performance, since the mean percentage estimation error in the worst
case (estimator (4) Model 2) is 1\% with low standard deviation.





\section{Conclusions} In this paper we introduce a new estimator of
 the diffusion part $IC$ and of the co-jumps in the quadratic
covariation of two semimartingales $X^{(q)}$.
To capture the separate contributions to the quadratic covariation  has important
applications in finance (forecasting, option pricing, risk and credit risk management).\\
The estimator $IC_{T,n}$ is constructed using a threshold
criterion introduced in \cite{Man01}, and consists in summing
properly selected cross products of increments of the two
processes. Our estimator is consistent, and when the two  jump
parts
 have only finite activity a joint CLT for $\hat{IC}_{T,n}$ and the estimators $\hat{IV}_{T}^{(q)}$ of the integrated variances
 is proved and delivers the following important consequences.\\
1. $\hat{IC}_{T,n}$ is also asymptotically Gaussian with  speed of
convergence
 $\sqrt{h}$. A central limit theorem in presence of
infinite activity jump parts is studied in a further paper
(\cite{GobMan07}) where we find that the speed of convergence of
$\hat{IC}_{T,n}$ is determined both by the dependence structure
between
the two processes $X^{(q)}$ and by the amount of jump activity of each $J^{(q)}$.\\
2. Consistent estimators both of the sum of the co-jumps
occurred within $[0,T]$ and  of each single co-jump are obtained.\\
3. We construct asymptotically Gaussian estimators of the regression coefficients $\beta$s and
of the correlation coefficient between the two processes $X^{(q)}$.\\
Further we find that in presence of FA jumps a slight modification
of $\hat{IC}_{T,n}$ is consistent even
when only non-synchronous observations are available.\\
We assess the choice of the threshold and check the performance of
our estimators on two different kind of simulated models which are
common in the financial literature. Model 1 has components with
stochastic volatilities and FA jumps, while Model 2 has components
with constant volatilities and IA jumps. We find that even with
five minutes observations the performances of the estimators of
$\sum_{0\leq t\leq T}\Delta J^{(1)}_{t}\Delta J^{(2)}_{t}$ and of
the single co-jumps are satisfactory. $\hat{IC}_{T, n}$ is satisfactory
in Model 1 while is biased in Model 2 but the corresponding
normalized bias has still Gaussian behavior.
\\

{\bf Acknowledgements.} The Authors thank Rama Cont for the important initial input to begin the work, Roberto Renò
and an anonymous referee for the important comments that helped to improve
and deepen our analysis.\\
The authors benefited of financial support by
Italian Goverment, grant MIUR 2006 n.206132713-001.



\section{Appendix}

 The following theorem is the key result, in the finite jump activity case,
 validating the idea that if $\DXq$ is larger than $\rh$ then some jumps occurred in $\Ii$ (and vice-versa). It is stated
 in the general case of not equally spaced observations.

\begin{theorem}\label{IndicatorNoJumps} (\cite{mancini07}, FA jumps)
Under the assumptions from
\textbf{A1} to  \textbf{A4}, given
 an arbitrary  partition
$\{t_0=0, t_1, ..., t_n=T\}$ of $[0,T],$ then for sufficiently small, but strictly positive,
$h:=\sup_{j=1..n} |t_{j}-t_{j-1}|$ (depending on $\omega$) we have  a.s.
$$ 1_{\{(\Delta_{j}X^{(q)})^{2}\leq
r_{h}\}}=1_{\{\Delta_{j}N^{(q)}=0\}},\quad j=1,2,....,n,\quad q=1,2.
$$

\vspace{-1.2cm}\qed\\
\end{theorem}

\textbf{\emph{Proof of Theorem \ref{JointCv} [Joint CLT] }} By
Theorem \ref{IndicatorNoJumps} we have, for all $t\in[0,T]$,
$$h^{-1/2}\left[ \sum_{t_j\leq t} \DXa_\star\DXb_\star- IC_t\right]\!=\!
h^{-1/2}\left[ \sum_{t_j\leq t} \DXa \DXb I_{\{ \DNa=0, \DNb=0\}}-IC_t\!\right] $$
$$
=h^{-1/2}\left[\sum_{t_j\leq t} \Delta_{j}D^{(1)}\Delta_{j}D^{(2)}-
 \int_{0}^{t}\rho_{s}\sigma_{s}^{(1)}\sigma_{s}^{(2)}ds\right]
-h^{-1/2}\sum_{t_j\leq t}
  \Delta_{j}D^{(1)}\Delta_{j}D^{(2)}1_{\{\Delta_{j}N^{(2)} \neq 0\}}
$$
$$
-h^{-1/2}\sum_{t_j\leq t}
\Delta_{j}D^{(1)}1_{\{\Delta_{j}N^{(1)}\neq 0\}}\Delta_{j}D^{(2)}
+h^{-1/2}\sum_{t_j\leq t}\Delta_{j}D^{(1)}1_{\{\Delta_{j}N^{(1)}\neq
0\}}\Delta_{j}D^{(2)} 1_{\{\Delta_{j}N^{(2)}\neq 0\}}.
$$
Each one of last three sums tends a.s. to zero as $h\ri 0$, since it contains at least one $I_{\{ \DN\neq 0 \}}$ and for
 any $q=1,2$ we have
$$
\Plim\Big|h^{-1/2}\sum_{t_j\leq t}
  \Delta_{j}D^{(1)}\Delta_{j}D^{(2)}1_{\{\Delta_{j}N^{(q)} \neq 0\}}\Big|
\leq   \Plim
K_{1}(\omega)K_{2}(\omega)\sqrt{h}\log\frac{1}{h}N_{T}^{(q)}=0.
$$
Moreover, analogously as in \cite{mancini07},
 for each $q=1,2$ we reach that
$$h^{-1/2}\left( \sum_{t_j\leq t}\DXq_\star- {IV}^{(q)}_t\right)=
h^{-1/2}\left(  \sum_{t_j\leq t}\DXq I_{\{ \DN=0\}}- {IV}^{(q)}_t\right) =$$
$$ h^{-1/2}\left(  \sum_{t_j\leq t} (\DD)^2 - {IV}^{(q)}_t\right)- h^{-1/2}\sum_{t_j\leq t} (\DD)^2 I_{\{ \DN\neq 0\}},$$
where the last term  tends a.s. to zero as $h\ri 0.$ Therefore we
have  that
\begin{equation}\label{dimteoJointCLTmatErrStima}
 h^{-1/2}\left(\begin{array}{cc}
\hat{IV}^{(1)}_n- {IV}^{(1)} & \hat{IC}_{T,n}- IC\\
\hat{IC}_n -IC & \hat{IV}_n^{(2)}- {IV}^{(2)}
\end{array}\right)
\eeq
has the same limit in distribution as
$$h^{-1/2}\left(\begin{array}{cc}
\sum_{t_j\leq t} (\DDa)^2- {IV}^{(1)} &  \sum_{t_j\leq t}\DDa\DDb- IC\\
\sum_{t_j\leq t} \DDa\DDb -IC & \sum_{t_j\leq t} (\DDb)^2- {IV}^{(2)}
\end{array}\right).$$
Note that
\beq\label{ItoFormperQCov}
\sum_{t_j\leq t} \DDa\DDb - IC_t= \sum_{t_j\leq t} \Big(\DDa\DDb - \Delta_j<\Da, \Db>\Big)\eeq
and, along the lines of
\cite{BS07Variation} (proof of Theorem 1, sec. 3.1),  using Itô formula we know that
$$ d(\Da\Db)= \Da_- d \Db + \Db_- d \Da + d<\Da,\Db>,$$
so $$ \Delta_j (\Da\Db) = \intIi\Da_{s-} d \Db_s + \intIi\Db_{s-}
d \Da_s + \Delta_j\!\!<\Da,\Db>.$$ Therefore $$\DDa\DDb= \Delta_j
(\Da\Db) -  \Da_{t_{j-1}}\DDb- \Db_{t_{j-1}}\DDa $$
$$=\intIi\Da_{s-} d \Db_s + \intIi\Db_{s-} d \Da_s + \Delta_j\!\!<\Da,\Db>-   \Da_{t_{j-1}}\DDb- \Db_{t_{j-1}}\DDa,$$
so that (\ref{ItoFormperQCov}) equals
$$ \int_0^t \Big(\Da_{s-} - \sum_{t_j\leq t} \Da_{t_{j-1}}I_{\{s\in \Ii\}} (s) \Big)d\Db_s +
\int_0^t \Big(\Db_{s-} - \sum_{t_j\leq t} \Db_{t_{j-1}}I_{\{s\in \Ii\}} (s) \Big)d\Da_s$$
$$= A^{(n)}_{12, t} + A^{(n)}_{21,t}, $$
where
$$ A^{(n)}=\left(\begin{array}{cc}
\int_0^\cdot \Big(\Da_{s-} - \Da_{\frac{[ns-]}{n}}\Big) d\Da_s &
      \int_0^\cdot \Big( \Da_{s-} - \Da_{\frac{[ns-]}{n}}\Big)d\Db_s \\
 \int_0^\cdot \Big(\Db_{s-} - \Db_{\frac{[ns-]}{n}}\Big)d\Da_s
      & \int_0^\cdot\Big( \Db_{s-} - \Db_{\frac{[ns-]}{n}} \Big)d\Db_s
\end{array}\right).$$
As special cases, for each $q=1,2$
$$ \sum_{t_j\leq t} \DDq - {IV}^{(q)}_t=  \sum_{t_j\leq t} \Big(\DDq - \Delta_j <\D,\D>\Big)=$$
$$ 2\int_0^t \Big(\D_{s-} - \sum_{j=1}^n \D_{t_{j-1}}I_{\{s\in\Ii\}} (s) \Big)d\D_s = 2A^{(n)}_{qq, t}.$$
By Theorem 5.5 in \cite{JacPro98} we have that
$$h^{-1/2} A^{(n)} \stackrel{st}{\ri} \frac{{\bf Z}}{\sqrt 2},$$
with ${\bf Z}$ as in (\ref{AsLawnelJointCLT}). It follows that, as $n\ri \infty$, (\ref{dimteoJointCLTmatErrStima})
converges stably in law to
$$ \frac{1}{\sqrt 2}\left(\begin{array}{cc}
2Z_{11} & Z_{12} + Z_{21}\\
Z_{12} + Z_{21} & 2Z_{22}
\end{array}\right).$$

 \vspace{-1cm}\qed

\textbf{\emph{Proof of Corollary \ref{speedCvBetasRho} [Speed of
convergence of $\beta$s and $\rho$, FA jumps]}}
For all $t\in[0,T]$ we have
$$h^{-1/2}\left(\hat{\beta}^{(i,j)}_{t,n} - \frac{IC_t}{{IV}^{(j)}_t }\right)=
h^{-1/2}\frac{\hat{IC}_{t,n} - IC_t}{\hat{IV}_{t,n}^{(j)}}+ h^{-1/2} IC_t \frac{{IV}^{(j)}_t-
\hat{IV}^{(j)}_{t,n}}{\hat{IV}^{(j)}_{t,n} {IV}^{(j)}_t},$$
therefore
$$h^{-1/2}\left(\hat{\beta}^{(i,j)}_{n} - \frac{IC}{{IV}^{(j)} }\right) \stackrel{st}{\ri} \frac{Z_{12}+ Z_{21}}{\sqrt 2 \ {IV}^{(j)}} - IC \frac{\sqrt 2 Z_{jj}}{ \ ({IV}^{(j)})^2}.$$
As for $\hat{\rho}^{(1,2)}_n$, note preliminarily that Theorem \ref{JointCv}
implies that $h^{-1/2}\left(\sqrt{{IV}^{(j)}} - \sqrt{\hat{IV}^{(j)}_n} \right)$ converges stably,
since, $t$ by $t$, $$h^{-1/2}\left(\sqrt{{IV}^{(j)}} - \sqrt{\hat{IV}^{(j)}_n} \right)=
\frac{h^{-1/2}\left({IV}^{(j)} - \hat{IV}^{(j)}_n \right)}{\sqrt{{IV}^{(j)}} + \sqrt{\hat{IV}^{(j)}_n}}
 \stackrel{st}{\ri} -\frac{Z_{jj}}{\sqrt{2 \ {IV}^{(j)}}}.$$
As a  consequence
$$h^{-1/2}\left(\hat{\rho}^{(1,2)}_n - \frac{IC}{\sqrt{{IV}^{(1)}{IV}^{(2)}}}\right)=$$
$$h^{-1/2}\frac{\hat{IC}_n - IC}{\sqrt{\hat{IV}^{(1)}_n\hat{IV}^{(2)}_n}}+
 h^{-1/2} IC \left(\frac{1}{\sqrt{\hat{IV}^{(1)}_n\hat{IV}^{(2)}_n}} - \frac{1}{\sqrt{{IV}^{(1)}{IV}^{(2)}}}\right).$$
The first term converges stably to $\frac{Z_{12}+ Z_{21}}{\sqrt{ 2 \ {IV}^{(1)}{IV}^{(2)}}}$, while the
second term equals
$$ \frac{h^{-1/2} IC }{\sqrt{\hat{IV}^{(1)}_n}}\left( \frac{1}{\sqrt{\hat{IV}^{(2)}_n}}-
\frac{\sqrt{\hat{IV}^{(1)}_n}}{\sqrt{{IV}^{(1)}}} \frac{1}{\sqrt{{IV}^{(2)}}}\right) =$$
$$ \frac{h^{-1/2} IC }{\sqrt{\hat{IV}^{(1)}_n}}\left(\frac{1}{\sqrt{\hat{IV}^{(2)}_n}}-
\frac{1}{\sqrt{{IV}^{(2)}}}\right) +
 \frac{h^{-1/2} IC }{\sqrt{\hat{IV}^{(1)}_n {IV}^{(2)} }}\left(1-\frac{\sqrt{\hat{IV}^{(1)}_n}}{\sqrt{{IV}^{(1)}}} \right)=$$
$$\frac{h^{-1/2} IC }{\sqrt{\hat{IV}^{(1)}_n}}
\left(\frac{\sqrt{{IV}^{(2)}}-\sqrt{\hat{IV}^{(2)}_n}}{\sqrt{\hat{IV}^{(2)}_n {IV}^{(2)}}}\right)+
\frac{h^{-1/2} IC }{\sqrt{\hat{IV}^{(1)}_n {IV}^{(2)} {IV}^{(1)} }}\left(\sqrt{{IV}^{(1)}} -\sqrt{\hat{IV}^{(1)}_n }\right)
$$
$$\stackrel{st}{\ri}
-\frac{ Z_{22} IC}{\sqrt{2\  {IV}^{(1)}} ({IV}^{(2)})^{3/2} }-
\frac{Z_{11} IC }{\sqrt{2\  {IV}^{(2)}} ({IV}^{(1)})^{3/2} }.$$
\qed

\emph{\textbf{Proof of Proposition \ref{DenominPerCLTFA}}}
\emph{\textbf{ [Estimate of the standard error for $\hat{IC}_n$,
FA jumps]}} For $t=T$ it is sufficient to show that as $n\ri
\infty$
$$
\tilde{v}^{(n)}_{2,2}(X^{(1)},X^{(2)})_{T}\stackrel{P}{\longrightarrow}
\int_{0}^{T}(2\rho_{t}^{2}+1)(\sigma_{t}^{(1)})^{2}(\sigma_{t}^{(2)})^{2}dt,
$$
and
$$
\tilde{w}^{(n)}(X^{(1)},X^{(2)})_{T}\stackrel{P}{\longrightarrow}
\int_{0}^{T}\rho_{t}^{2}(\sigma_{t}^{(1)})^{2}(\sigma_{t}^{(2)})^{2}dt.
$$
For $t<T$ the proof is analogous with $\sum_{j=1}^{n}$ replaced by
$\sum_{j:t_j \leq t}.$
 By
Theorem \ref{IndicatorNoJumps} we can write
$$
  \Plim\ \tilde{v}^{(n)}_{2,2}(X^{(1)},X^{(2)})_{T}=
   \Plim\ h^{-1}\sum_{j=1}^{n}(\Delta_{j}D^{(1)})^2 1_{\{\Delta_{j}N^{(1)}=0\}}
   (\Delta_{j}D^{(2)})^2 1_{\{\Delta_{j}N^{(2)}=0\}}$$
  $$
=\Plim\ v^{(n)}_{2,2}(D^{(1)},D^{(2)})_{T}- \Plim\
h^{-1}\sum_{j=1}^{n}(\Delta_{j} D^{(1)})^{2}
   (\Delta_{j} D^{(2)})^{2}1_{\{\Delta_{j}N^{(1)}\neq
   0\}}
$$
$$
-\Plim\ h^{-1}\sum_{j=1}^{n}(\Delta_{j} D^{(1)})^{2}
   (\Delta_{j} D^{(2)})^{2}1_{\{\Delta_{j}N^{(2)}\neq 0\}}
$$
$$
+\Plim\ h^{-1}\sum_{j=1}^{n}(\Delta_{j} D^{(1)})^{2}
   1_{\{\Delta_{j}N^{(1)}\neq 0\}}
   (\Delta_{j} D^{(2)})^{2}1_{\{\Delta_{j}N^{(2)}\neq 0\}}.
$$
By Theorem 2.1 in \cite{BarGraJacPodShe05}, $$\Plim\
v^{(n)}_{2,2}(D^{(1)},D^{(2)})_{T}=\int_{0}^{T}(2\rho_t^{2}+1)
(\sigma_{t}^{(1)})^{2}(\sigma_{t}^{(2)})^{2}dt,$$ whereas the
other terms are all zero. In fact for any $q=1,2$
\beq\label{TermCoiSalTrasc}
  \Plim\ h^{-1}\sum_{j=1}^{n}(\Delta_{j} D^{(1)})^{2}
   (\Delta_{j}D^{(2)})^{2}1_{\{\Delta_{j}N^{(q)}\neq
   0\}}
\leq\Plim\
K_{1}^{2}(\omega)K_{2}^{2}(\omega)h\Big(\log\frac{1}{h}\Big)^{2}N_{T}^{(q)}=0.
\eeq
Now we deal with $\tilde{w}^{(n)}(X^{(1)},X^{(2)})_{T}$.
Analogously as before
$$
   \Plim\ \tilde{w}^{(n)}(X^{(1)},X^{(2)})_{T}
$$
$$
=\Plim\ h^{-1}\sum_{j=1}^{n-1}
  \left[\prod_{i=0}^{1}\Delta_{j+i}D^{(1)}(1-1_{\{(\Delta_{j+i}N^{(1)}\neq
  0\}})
  \prod_{i=0}^{1}\Delta_{j+i}D^{(2)}(1-1_{\{(\Delta_{j+i}N^{(2)}\neq
  0\}})\right],
$$
which coincides with the sum of $ \Plim\
w^{(n)}(D^{(1)},D^{(2)})_{T}$ with a finite number of  terms which
are shown to be negligible. By Theorem 2.1 in
\cite{BarGraJacPodShe05}, $ \Plim\ w^{(n)}(D^{(1)},D^{(2)})_{T}$ $
=\int_{0}^{T}\rho_{t}^{2}(\sigma_{t}^{(1)})^{2}(\sigma_{t}^{(2)})^{2}dt$,
while the  other terms are given by the product of
$\prod_{i=0}^{1}\Delta_{j+i}D^{(1)}$
$\prod_{i=0}^{1}\Delta_{j+i}D^{(2)}$ with at least one of the
indicators
 $1_{\{\Delta_{j+s}N^{(q)}\neq 0\}}$, for an $s\in\{0,1\}$. Therefore the limit in
probability of each such term is zero as in (\ref{TermCoiSalTrasc}). \qed


\vspace{1cm}
\emph{\textbf{Proof of Corollary \ref{CLTFiniteActivity}  [CLT for the standardized version of $\hat{IC}_{t,n} - IC_t$,
FA jumps]}}.
By Theorem \ref{JointCv} we have  $$
h^{-1/2}\left(\hat{IC}_n-IC\right) \stackrel{st}{\ri} $$
$$ \frac{1}{\sqrt{2}}\left(\int_0^\cdot 2\sigma_s^{(1)}\sigma_s^{(2)}\rho_s dB_s^{11}+
\int_0^\cdot \sigma_s^{(1)}\sigma_s^{(2)}\sqrt{1-\rho_s^2} \Big[dB_s^{12}+dB_s^{21}\Big]
\right).
$$
The variance of the last term at time $t$ is
$\int_{0}^{t}(1+\rho_{s}^{2})(\sigma_{t}^{(s)})^{2}(\sigma_{s}^{(2)})^{2}ds$.
By Proposition \ref{DenominPerCLTFA} we then obtain that
$$
\frac{\hat{IC}_{t,n}-IC_t}
{\sqrt{h}\sqrt{\tilde{v}^{(n)}_{2,2}(X^{(1)},X^{(2)})_{t}-\tilde{w}^{(n)}(X^{(1)},X^{(2)})_{t}}}\stackrel{d}{\ri}
{\cal N},
$$
where ${\cal N}$ is a standard Gaussian r.v.. \qed

\vspace{1cm}\emph{\textbf{Proof of Theorem \ref{asyncObs}
[Asynchronous observations]}}  Note that we can assume that $a$
and $\sigma$ are bounded on $[0,T]$ (\cite{Jac08}), so that the
Brownian semimartingale parts $D^{(q)}$ of $X^{(q)}$ belong to
$L^{8}$. Using Theorem \ref{IndicatorNoJumps} in the not
evenly-spaced observations case (\cite{mancini07}) with
$h:=\sup_{j=1..m^{(n)}} (\tau_{j}- \tau_{j-1}) \vee \sup_{i=
1..k^{(n)}} (\nu_{i}- \nu_{i-1})$, a.s. for sufficiently small $h$
we can write
$$\sum_{j=1..m^{(n)},\  i= 1..k^{(n)}}
(X^{(1)}_{\tau_{j}}-X^{(1)}_{\tau_{j-1}})_\star\
(X^{(2)}_{\nu_{i}}-X^{(2)}_{\nu_{i-1}})_\star\
1_{\{]\tau_{j-1}, \tau_{j}]\cap]\nu_{i-1}, \nu_{i} ] \neq
\emptyset \}}=$$
$$\!\! \!\!\!\!\sum_{j=1..m^{(n)},\  i= 1..k^{(n)}}\!\!\!\!
(D^{(1)}_{\tau_{j}}-D^{(1)}_{\tau_{j-1}})
1_{\{N^{(1)}_{\tau_{j}}-N^{(1)}_{\tau_{j-1}}=0\}}
(D^{(2)}_{\nu_{i}}-D^{(2)}_{\nu_{i-1}})
1_{\{N^{(2)}_{\nu_{i}}-N^{(2)}_{\nu_{i-1}}=0\}}
1_{\{]\tau_{j-1}, \tau_{j}]\cap]\nu_{i-1}, \nu_{i} ] \neq
\emptyset \}}$$
$$= \sum_{j=1..m^{(n)},\  i= 1..k^{(n)}}\!\!\!\!
(D^{(1)}_{\tau_{j}}-D^{(1)}_{\tau_{j-1}})
(D^{(2)}_{\nu_{i}}-D^{(2)}_{\nu_{i-1}}) 1_{\{]\tau_{j-1},
\tau_{j}]\cap]\nu_{i-1}, \nu_{i} ] \neq \emptyset \}}$$
$$ -\sum_{j=1..m^{(n)}, \ i= 1..k^{(n)}}\!\!\!\!
(D^{(1)}_{\tau_{j}}-D^{(1)}_{\tau_{j-1}})
(D^{(2)}_{\nu_{i}}-D^{(2)}_{\nu_{i-1}}) \left[1_{\{
N^{(1)}_{\tau_{j}}-N^{(1)}_{\tau_{j-1}}\neq 0\}} + 1_{\{
N^{(2)}_{\nu_{i}}-N^{(2)}_{\nu_{i-1}}\neq 0\}}+\right.$$
$$\left.1_{\{ N^{(1)}_{\tau_{j}}-N^{(1)}_{\tau_{j-1}}\neq 0, N^{(2)}_{\nu_{i}}-N^{(2)}_{\nu_{i-1}}\neq 0\}}
\right] 1_{\{]\tau_{j-1}, \tau_{j}]\cap]\nu_{i-1}, \nu_{i} ] \neq
\emptyset \}},$$ The first sum of the r.h.s. tends to $IC_T$ in
probability by Corollary 2.2 in (\cite{HayKus08}), with $f\equiv g
\equiv 1$, while each sum in the second term is dominated in
absolute value, for a suitable $q$, by
$$ \sup_{j}|D^{(1)}_{\tau_{j}}-D^{(1)}_{\tau_{j-1}}|\sup_i |D^{(2)}_{\nu_{i}}-D^{(2)}_{\nu_{i-1}}| N^{(q)}_T,
$$ which tends a.s. to zero as $h\ri 0$,  by Lemma \ref{ModContD}. \qed

\vspace{0.5cm}
 The following facts are used within the proof of Theorem
\ref{consistencyIA}.


Without loss of generality (as in \cite{Jac08}, Lemma 4.6) we can assume that
     \begin{center}
     \begin{minipage}[c]{0.9\linewidth}
     \textbf{A5}.
$\int_{x\in\R} 1\wedge (\gamma^{(q)})^2(\omega, t,x) dx$ is bounded.
 \end{minipage}
 \end{center}

\begin{lemma}\label{IDX0magredEDXdq} For each $q=1,2$ we have the following.
\begin{enumerate}
\item If processes $a$ and $\sigma$ are càdlàg then, under {\bf A3}, a.s., for small $h$,
    $1_{\{(\Delta_{j}D^{(q)})^{2}>r_{h}\}}=0$, uniformly in $j$.
\item Under \textbf{A5} we have that, for each $j=1,.., n$,
$E[\DXdq]\leq K h$, for a positive constant $K$.
\end{enumerate}
\end{lemma}

\dimo \ Part 1. is a consequence of Lemma \ref{ModContD}.\\
Part 2. $$E[\DXdq]= E[\intIi\int_{ |\gamma^{(q)}|\leq
1}(\gamma^{(q)})^2 \underline\nu^{(q)}(dx, ds) ]= E[\intIi\int_{
|\gamma^{(q)}|\leq 1}(\gamma^{(q)})^2 dx ds ]:$$ since, by
assumption {\bf A5}, $\int_{|\gamma^{(q)}|\leq 1} (\gamma^{(q)})^2
dx$ is bounded, the last term is dominated by $K h$ for some
positive constant $K$.\qed\\

The following lemma generalizes analogous results given in
\cite{mancini07} from the framework of Lévy jumps to the one of
Itô semimartingale jumps.
\begin{lemma}\label{4puntiPerDimoConsIA} The following facts hold.
\begin{enumerate}
\item Let us consider any
sequence $\pi_n$ of partitions $\{ 0, t_1,.., t_n=T\}$ of $[0,T]$,
$n\in\N$,  such that $\max_{j=1..n} |t_j - t_{j-1}| \ri 0$ as $n\ri
\infty$. For each $q=1,2,$ as long as $\tilde J_2^{(q)}$ is a
semimartingale, we can find a subsequence $n_k$ for which a.s., for
any $\delta >0$ there exists a sufficiently large $k$ such that for
all $j=1, .., n_k$ on $\{\DXdq\leq 4 r_{h_k}\}$ we have
$$(\Delta \tilde J_{2,s}^{(q)})^2\leq 4r_{h_k}+\delta, \quad \forall s\in \Ii.$$
\item Under {\bf A3} and {\bf A5}, for each $q=1,2$, we have
    $\sum_{j=1}^n P\{ \DN\neq 0, \DXdq>4\rh \}\ri 0 $ as $h\ri 0$.
\end{enumerate}
\end{lemma}

\textbf{\emph{Proof}}.  {\it Statement 1}  is a consequence of the
fact that (\cite{Met82}, Theorem 25.1) there is a subsequence $n_k$
such that, defined $h_k=T/n_k$, $\sum_{j=1}^{[t/h_k]} \DXdq$ tends
to $\sum_{s\in [0,t]} (\Delta\tilde J_{2, s}^{(q)})^2$ a.s.
uniformly w.r.t. $t\in [0,T]$, as $k\ri \infty$, where $[x]$ denotes
the integer part of $x$. Since a.s.
$$\sup_{j=1..n_k}|\DXdq-\sum_{s\in \Ii} (\Delta\tilde J_{2, s}^{(q)})^2|$$ $$=
\sup_{j=1..n_k} \left\{\big[\sum_{\ell=1}^{[t_j/h_k]} (\Delta_\ell
\tilde J_2)^2-\sum_{s\in [0,t_j]} (\Delta\tilde J_{2,s}^{(q)})^2\big]-
\big[ \sum_{\ell=1}^{[t_{j-1}/h_k]} (\Delta_\ell
\tilde J_2)^2-\sum_{s\in [0,t_{j-1}]} (\Delta\tilde J_{2,
s}^{(q)})^2\big]\right\}
$$ $$\leq 2 \sup_{t\in [0,T]} \big[\sum_{j=1}^{[t/h_k]}
\DXdq-\sum_{s\in [0,t]} (\Delta\tilde J_{2, s}^{(q)})^2\big]\ri 0,
$$ we in fact have that a.s. for all $j=1, .., n_k$ each
squared increment $\DXdq$ is uniformly, on $j$, arbitrarily close to
$\sum_{s\in \Ii} (\Delta\tilde J_{2, s}^{(q)})^2$. More precisely,
a.s.  for all $\delta>0$ we can find a sufficiently large $k$ such
that
$$\sup_{j=1..n_k}\left|\DXdq-\sum_{s\in \Ii} (\Delta\tilde J_{2, s}^{(q)})^2\right|<\delta, $$
so, for all $j$ such that $\DXdq\leq 4\rh$ we have $$ \sum_{s\in \Ii} (\Delta\tilde
J_{2, s}^{(q)})^2\leq \sup_{j=1..n_k}\left|\DXdq-\sum_{s\in \Ii}
(\Delta\tilde J_{2, s}^{(q)})^2\right|+ \DXdq \leq 4\rh+\delta.$$
In particular for any $s\in \Ii$ with $j$ such that $\DXdq\leq 4\rh$, each squared jump size
$(\Delta\tilde J_{2, s}^{(q)})^2$
is bounded by $4\rh+\delta$.\\
{\it Statement  2}. The predictable
compensator of $N^{(q)}_t=\sum_{s\leq t} I_{|\Delta J^{(q)}_s|>1}$ and the predictable
quadratic variation of $\tilde J_2^{(q)}$ are of the form $\Lambda^{(q)}_t= \int_0^t
\lambda^{(q)}_s ds$ and $\Lambda'^{(q)}_t= \int_0^t \lambda'^{(q)}_s ds$ respectively.
Assumption {\bf A5} guarantees that both
$\int_{|\gamma^{(q)}|>1} 1\  dx$ and \\ $\int_{|\gamma^{(q)}|\leq 1} (\gamma^{(q)})^2\  dx$
are bounded, and therefore that $\lambda^{(q)}$ and $\lambda'^{(q)}$ are bounded processes.
Using exactly the same argument as in
\cite{AitJacStimBet07}, eq. (60), with $\delta=1$ and $\zeta=3\sqrt{r_h}$ and replacing
$M(\delta)=\int \int_{|x|\leq 1} x \tilde \mu (dx, dt)$ with our
 $\tilde J_2^{(q)}$, we conclude that
$$\sum_{j=1}^nP\{ \DN\neq 0, \DXdq>4\rh \}=
O(nh \ \frac{h}{\rh}). \qed$$

\vspace{1cm} For any $\delta >0$ denote by $Z_{h_k}^{(q), 
\delta}$ the following pure jump plus drift semimartingales having
only jumps bounded in absolute value by $\sqrt{4r_{h_k}+\delta}$, $q=1,2$:
$$
Z_{h_k,t}^{(q), \delta}:=
\int_0^t\int_{|\gamma^{(q)}|\leq \sqrt{4r_{h_k}+\delta }}\gamma^{(q)} \underline{\tilde\mu}^{(q)}(dx,ds)
-\int_0^t\int_{\sqrt{4r_{h_k}+\delta}< |\gamma^{(q)}| \leq 1}\ \gamma^{(q)}\ dx dt,\quad t\geq 0.
$$
By Lemma  \ref{4puntiPerDimoConsIA} we have that for any
$\delta>0,$ for sufficiently large $k$  the indexes $j$ for which
$\DXdq\leq 4r_{h_k}$ are such that the increment $\DXdq$ coincides
with the increment $(\Delta_jZ_{h_k}^{(q), \delta})^2$ of
$Z_{h_k}^{(q), \delta}$, since $\DXd$ does not contain jumps
larger than $\sqrt{4r_{h_k}+\delta}$.

\begin{lemma}\label{SumDXddoveleqrrizero} 
For each $q=1,2$
$$\Plim \sum_{j=1}^{n}
(\Delta_{j}\tilde{J}^{(q)}_{2})^21_{\{(\Delta_{j}\tilde{J}^{(q)}_{2})^{2}\leq
4r_{h}\}}=0$$
\end{lemma}

\emph{\textbf{Proof}} Consider the sequence of partitions
$\pi_n=\{0, T/n, 2 T/n, .., T \}$. Take any subsequence
$\pi_{n_\ell}$. By Lemma \ref{4puntiPerDimoConsIA}, point {\it 1},
a.s. there exists a sub-subsequence $n_{\ell_k}$ such that for any
$\delta>0$ and $k$ sufficiently large then for all $j=1, ..,
n_{\ell_k}$ on $\DXdq\leq 4 r_{h_{\ell_k}}$ we have
$\DXd=(\Delta_jZ_{h_{\ell_k}}^{(q), \delta})^2$. Denote
$$S^{(q)}_n:= \sum_{j=1}^{n}
(\Delta_{j}\tilde{J}^{(q)}_{2})^21_{\{(\Delta_{j}\tilde{J}^{(q)}_{2})^{2}\leq
4r_{h}\}}.$$ Therefore
$$0\leq \Plimk S^{(q)}_{n_{\ell_k}}
\leq \Plimk \sum_{j=1}^{n_{\ell_k}}(\Delta_j Z_{h_{\ell_k}}^{(q),
\delta})^2$$
$$= \Plimk \int_0^T\int_{|\gamma^{(q)}|\leq \sqrt{4r_{h_{\ell_k}}+\delta }} \ (\gamma^{(q)})^2 \underline{\nu}^{(q)} (dx,
ds)=\int_0^T\int_{|\gamma^{(q)}|\leq \sqrt{\delta }} \ (\gamma^{(q)})^2 dx ds.$$
Since a.s. $\int_{|\gamma^{(q)}|\leq 1} \ (\gamma^{(q)})^2
dx<\infty$, the last term above
tends a.s. to zero as $\delta \ri 0$, which implies that  $\Plimk S^{(q)}_{n_{\ell_k}}=0$.\\
Since then from any subsequence of $S_n$ we can extract a
sub-subsequence tending to zero in probability, we in fact have
that the whole sequence $S^{(q)}_n\ri0$ in probability, as we need.\qed


\vspace{0.5cm}
\emph{\textbf{Proof of Theorem \ref{consistencyIA}}}. 
 We decompose
$\hat{IC}_{T,n}-IC_T$ into the sum of five terms and
we show that each term tends a.s. to zero, as $n\rightarrow
\infty$.  We need some further notation. Recall that for each $q=1,2$
$$ D^{(q)}_{t} = \int_{0}^{t}a^{(q)}_{s}ds+\int_{0}^{t}\sigma^{(q)}_{s}dW_{s}^{(q)},  $$
and denote
$$
   Y^{(q)}_{t} := D^{(q)}_{t}+J_{1t}^{(q)},
$$
so that we have $X^{(q)}_{t}=Y^{(q)}_{t}+\tilde{J}_{2t}^{(q)}$, $q=1,2$.

Adding and subtracting $\sum_{j=1}^{n}
  \Delta_{j}Y^{(1)}1_{\{(\Delta_{j}Y^{(1)})^{2}\leq
9r_{h}\}}\Delta_{j}Y^{(2)}1_{\{(\Delta_{j}Y^{(2)})^{2}\leq
9r_{h}\}}$ from $\hat{IC}_{T,n}-IC_T$,   we reach 
$$|\hat{IC}_{T,n}-IC_T|$$
$$ =\Big|\!\sum_{j=1}^{n}
  (\Delta_{j}Y^{(1)}\!+\! \Delta_{j}\tilde{J}_{2}^{(1)})\!1_{\{(\Delta_{j}X^{(1)})^{2}\leq r_{h}\}}
  (\Delta_{j}Y^{(2)}\!+\! \Delta_{j}\tilde{J}_{2}^{(2)})\!1_{\{(\Delta_{j}X^{(2)})^{2}\leq r_{h}\}}
 \!-IC_T \Big|
$$
$$
  \leq\Big|\sum_{j=1}^{n}
  \Delta_{j}Y^{(1)}1_{\{(\Delta_{j}Y^{(1)})^{2}\leq
9r_{h}\}}\Delta_{j}Y^{(2)}1_{\{(\Delta_{j}Y^{(2)})^{2}\leq
9r_{h}\}}-IC_T \Big|
$$
\begin{equation}\label{ConsIAtildev11menorhos1s2}
\begin{array}{c}
 \!\!\!\!\!+\Big|\sum_{j=1}^{n}\Delta_{j}Y^{(1)}\Delta_{j}Y^{(2)}\Big(1_{\{(\Delta_{j}X^{(1)})^{2}\leq
r_{h}\}}1_{\{(\Delta_{j}X^{(2)})^{2}\leq
r_{h}\}}\!-\!1_{\{(\Delta_{j}Y^{(1)})^{2}\leq
9r_{h}\}}1_{\{(\Delta_{j}Y^{(2)})^{2}\leq 9r_{h}\}}\Big)\Big|
\\
\\
+\Big|\sum_{j=1}^{n}\Delta_{j}Y^{(1)}\Delta_{j}\tilde{J}_{2}^{(2)}1_{\{(\Delta_{j}X^{(1)})^{2}\leq
r_{h}\}}1_{\{(\Delta_{j}X^{(2)})^{2}\leq r_{h}\}}\Big|
\\
\\
+\Big|\sum_{j=1}^{n}\Delta_{j}\tilde{J}_{2}^{(1)}\Delta_{j}Y^{(2)}1_{\{(\Delta_{j}X^{(1)})^{2}\leq
r_{h}\}}1_{\{(\Delta_{j}X^{(2)})^{2}\leq r_{h}\}}\Big|+
\end{array}
\end{equation}
$$+\Big|\sum_{j=1}^{n}\Delta_{j}\tilde{J}_{2}^{(1)}\Delta_{j}\tilde{J}_{2}^{(2)}1_{\{(\Delta_{j}X^{(1)})^{2}\leq
r_{h}\}}1_{\{(\Delta_{j}X^{(2)})^{2}\leq r_{h}\}}\Big|.
$$

\vspace{0.5cm}
\n The first term tends to zero in probability by Corollary
\ref{hatICconsistFA}. The second term 
coincides with
\begin{equation}\label{somma}
\begin{array}{c}
\Big|\sum_{j=1}^{n}\Delta_{j}Y^{(1)}\Delta_{j}Y^{(2)}\Big[1_{\{(\Delta_{j}X^{(1)})^{2}\leq
r_{h},(\Delta_{j}X^{(2)})^{2}\leq r_{h},(\Delta_{j}Y^{(1)})^{2}>
9r_{h}\}}
\\
\\
+1_{\{(\Delta_{j}X^{(1)})^{2}\leq
r_{h},(\Delta_{j}X^{(2)})^{2}\leq r_{h},(\Delta_{j}Y^{(2)})^{2}>
9r_{h}\}}
\\
\\
-1_{\{(\Delta_{j}X^{(1)})^{2}\leq
r_{h},(\Delta_{j}X^{(2)})^{2}\leq r_{h},(\Delta_{j}Y^{(1)})^{2}>
9r_{h},(\Delta_{j}Y^{(2)})^{2}> 9r_{h}\}}
\\
\\
\!\!\!\!-1_{\{(\Delta_{j}X^{(1)})^{2}>
r_{h},(\Delta_{j}Y^{(1)})^{2}\leq
9r_{h},(\Delta_{j}Y^{(2)})^{2}\leq
9r_{h}\}}-1_{\{(\Delta_{j}X^{(2)})^{2}>
r_{h},(\Delta_{j}Y^{(1)})^{2}\leq
9r_{h},(\Delta_{j}Y^{(2)})^{2}\leq 9r_{h}\}}
\\
\\
+1_{\{(\Delta_{j}X^{(1)})^{2}> r_{h},(\Delta_{j}X^{(2)})^{2}>
r_{h},(\Delta_{j}Y^{(1)})^{2}\leq
9r_{h},(\Delta_{j}Y^{(2)})^{2}\leq 9r_{h}\}}\Big]\Big|.\\
\end{array}
\end{equation}

\vspace{0.5cm} \n All these terms tend a.s. to zero. In fact for
the first three ones notice that on
$\{(\Delta_{j}X^{(q)})^{2}\leq r_{h},(\Delta_{j}Y^{(q)})^{2}>
9r_{h}\}$ we have $\sqrt{r_{h}}\geq
|\Delta_{j}X^{(q)}|\geq|\Delta_{j}Y^{(q)}|-|\Delta_{j}\tilde{J}^{(q)}_{2}|$
and thus $|\Delta_{j}\tilde{J}^{(q)}_{2}|\geq
|\Delta_{j}Y^{(q)}|-\sqrt{r_{h}}>3\sqrt{r_{h}}-\sqrt{r_{h}}=2\sqrt{r_{h}}$,
so \ that \ $\{(\Delta_{j}X^{(q)})^{2}$ $\leq
r_{h},(\Delta_{j}Y^{(q)})^{2}> 9r_{h}\}\subset \{(\Delta_{j}X^{(q)})^{2}\leq
r_{h},\DXdq>4r_{h}\}$,
$q=1,2$, and thus the probability that each one of the first three
terms of (\ref{somma}) is non zero is bounded by \beq\label{P}
P\Big\{\sum_{j=1}^{n} 1_{\{(\Delta_{j}X^{(q)})^{2}\leq
r_{h},\DXdq> 4r_{h}\}}\neq 0\Big\}\eeq
 for a
suitable $q$. Now on ${\{(\Delta_{j}X^{(q)})^{2}\leq r_{h},\DXdq >
4r_{h}\}}$ we in fact have that $\DN\neq 0$. Actually, since
$$ 2\sqrh -|\DY|<|\DXdq|- |\DY|\leq |\Delta_{j}X^{(q)}|\leq \sqrh$$
then
$$   K_q \sqrt{\mc} +|\DJu|\geq |\DD|+|\DJu|\geq |\DY|> \sqrh, $$
so $$|\DJu|>  \sqrh\left(1-K_q\sqrt{\frac{\mc}{\rh}}\right):$$ since a.s. for sufficiently small
$h$ the quantity $1-K_q \sqrt{\mc/\rh}$ is positive,
then in fact $|\DJu|>0$, so
that $\DN\neq 0.$\\
 So (\ref{P}) is dominated by
$\sum_{j=1}^n P\{\DN\neq 0,\DXdq> 4r_{h}\}$ which
tends to zero as $h\ri 0$ by Lemma \ref{4puntiPerDimoConsIA}, part 2).\\
As for the last three terms of (\ref{somma}) note that on
$\{(\Delta_{j}Y^{(q)})^{2}\leq 9r_{h}\}$ we have a.s., for $h$
small such that $\DN\in\{0,1\}$,
$$
\Delta_{j}N^{(q)}\leq
|\Delta_{j}J_{1}^{(q)}|= |\Delta_{j}Y^{(q)}
-\Delta_{j}D^{(q)}|<|\Delta_{j}D^{(q)}|+|\Delta_{j}Y^{(q)}|
$$
$$
\leq|\Delta_{j}D^{(q)}|+3\sqrt{r_{h}}\leq
K_q\sqrt{\mc}+3\sqrt{r_{h}}\rightarrow 0,\quad q=1,2,
$$
hence, for small $h$ on $\{(\Delta_{j}Y^{(q)})^{2}\leq 9r_{h}\}$
we have $\Delta_{j}N^{(q)}=0$, $j=1,..,n$. Therefore
$\{(\Delta_{j}X^{(q)})^{2}>r_{h},(\Delta_{j}Y^{(q)})^{2}\leq
9r_{h}\}\subset
\{(\Delta_{j}D^{(q)}+\Delta_{j}\tilde{J}_{2}^{(q)})^{2}>r_{h}\}\subset
\{(\Delta_{j}D^{(q)})^{2}>\frac{r_{h}}{4}\}\cup
\{(\Delta_{j}\tilde{J}_{2}^{(q)})^{2}>\frac{r_{h}}{4}\}$, $q=1,2$;
however, by Lemma \ref{IDX0magredEDXdq} part 1), a.s., for
small $h$, $ 1_{\{(\Delta_{j}D^{(q)})^{2}>\frac{r_{h}}{4}\}}=0, $
thus the last three terms of (\ref{somma}) are dominated by
$$\sum_{j=1}^{n}|\Delta_{j}D^{(1)}\Delta_{j}D^{(2)}|1_{\{\DXdq>\frac{\rh}{4}\}}$$
for a suitable $q$. However this last term tends to zero in
probability, since
$$\sum_{j=1}^{n}|\Delta_{j}D^{(1)}\Delta_{j}D^{(2)}|1_{\{\DXdq>\frac{\rh}{4}\}}\leq
K_1K_2\mc \sum_{j=1}^{n}1_{\{\DXdq>\frac{\rh}{4}\}}$$ and
$E[\mc \sum_{j=1}^{n}1_{\{\DXdq>\frac{\rh}{4}\}}]= \mc
\sum_{j=1}^n P\{\DXdq>\frac{\rh}{4}\}=O(\frac{\mc}{\rh}nh).$

 We now show that the third and fourth terms of the right hand
 side of
(\ref{ConsIAtildev11menorhos1s2}), which are similar, tend to zero
in probability. We have

\begin{equation}\label{sommaDeltaYDeltatildeJ2}
\begin{array}{c}
\sum_{j=1}^{n}\Delta_{j}Y^{(1)}\Delta_{j}\tilde{J}_{2}^{(2)}1_{\{(\Delta_{j}X^{(1)})^{2}\leq
r_{h}\}}1_{\{(\Delta_{j}X^{(2)})^{2}\leq r_{h}\}}
\\
\\
=\sum_{j=1}^{n}\Delta_{j}Y^{(1)}\Delta_{j}\tilde{J}_{2}^{(2)}\Big[1_{\{|\Delta_{j}X^{(1)}|\leq
\sqrt{r_{h}},|\Delta_{j}\tilde{J}_{2}^{(1)}|\leq 2\sqrt{r_{h}}\}}
1_{\{|\Delta_{j}X^{(2)}|\leq
\sqrt{r_{h}},|\Delta_{j}\tilde{J}_{2}^{(2)}|\leq 2\sqrt{r_{h}}\}}
\\
\\
+1_{\{|\Delta_{j}X^{(1)}|\leq
\sqrt{r_{h}},|\Delta_{j}\tilde{J}_{2}^{(1)}|> 2\sqrt{r_{h}}\}}
1_{\{|\Delta_{j}X^{(2)}|\leq
\sqrt{r_{h}},|\Delta_{j}\tilde{J}_{2}^{(2)}|> 2\sqrt{r_{h}}\}}
\\
\\
+1_{\{|\Delta_{j}X^{(1)}|\leq
\sqrt{r_{h}},|\Delta_{j}\tilde{J}_{2}^{(1)}|> 2\sqrt{r_{h}}\}}
1_{\{|\Delta_{j}X^{(2)}|\leq
\sqrt{r_{h}},|\Delta_{j}\tilde{J}_{2}^{(2)}|\leq 2\sqrt{r_{h}}\}}
\\
\\
+1_{\{|\Delta_{j}X^{(1)}|\leq
\sqrt{r_{h}},|\Delta_{j}\tilde{J}_{2}^{(1)}|\leq 2\sqrt{r_{h}}\}}
1_{\{|\Delta_{j}X^{(2)}|\leq
\sqrt{r_{h}},|\Delta_{j}\tilde{J}_{2}^{(2)}|>
2\sqrt{r_{h}}\}}\Big].
\\
\end{array}
\end{equation}

\n As before  on ${\{(\Delta_{j}X^{(q)})^{2}\leq r_{h},\DXdq >
4r_{h}\}}$ we have that $\DN\neq 0$, so, for each one of the last
three terms of (\ref{sommaDeltaYDeltatildeJ2}), the probability it
is different from zero is dominated by
$$\sum_{j=1}^n P\{ \DN\neq 0, \DXdq> 4\rh\}\ri 0.$$
Now we show that the first term of (\ref{sommaDeltaYDeltatildeJ2})
is asymptotically negligible. Notice that on
$\{|\Delta_{j}X^{(q)}|\leq
\sqrt{r_{h}},|\Delta_{j}\tilde{J}_{2}^{(q)}|\leq 2\sqrt{r_{h}}\}$
a.s. for small $h$ we have $\Delta N_j^{(q)}=0$; in fact a.s., for
small $h$ we have $\DN\in\{0,1\},$ and
$$
\Delta_{j}N^{(q)}\leq
|\Delta_{j}J_{1}^{(q)}|=|\Delta_{j}X^{(q)}-\Delta_{j}D^{(q)}-\Delta_{j}\tilde{J}^{(q)}_{2}|
$$
$$
\leq\sqrt{r_{h}}+\sup_{j}|\Delta_{j}D^{(q)}|+2\sqrt{r_h}\rightarrow
0,
$$
for all $j=1, .., n$, for $q=1,2$.
So we have
$$
\Plim\
\left|\sum_{j=1}^{n}\Delta_{j}Y^{(1)}\Delta_{j}\tilde{J}_{2}^{(2)}1_{\{|\Delta_{j}X^{(1)}|\leq
\sqrt{r_{h}},|\Delta_{j}\tilde{J}_{2}^{(1)}|\leq 2\sqrt{r_{h}}\}}
1_{\{|\Delta_{j}X^{(2)}|\leq
\sqrt{r_{h}},|\Delta_{j}\tilde{J}_{2}^{(2)}|\leq
2\sqrt{r_{h}}\}}\right|
$$
$$
\leq\Plim\sum_{j=1}^{n}|\Delta_{j}D^{(1)}\Delta_{j}\tilde{J}_{2}^{(2)}|
1_{\{|\Delta_{j}\tilde{J}_{2}^{(2)}|\leq 2\sqrt{r_{h}}\}}.
$$
By the Cauchy-Schwarz inequality, 
last term is dominated by
$$
\Plim \sqrt{\sum_{j=1}^{n}(\Delta_{j}D^{(1)})^2}
\sqrt{\sum_{j=1}^{n}(\Delta_{j}\tilde{J}_{2}^{(2)})^2
1_{\{|\Delta_{j}\tilde{J}_{2}^{(2)}|\leq 2\sqrt{r_{h}}\}}}
$$
$$
\leq\sqrt{\int_0^T (\sigma^{(1)}_s)^2 ds  }\ \Plim
\sqrt{S^{(2)}_n}=0,
$$
by Lemma \ref{SumDXddoveleqrrizero}.

It remains to consider the last term of
(\ref{ConsIAtildev11menorhos1s2}), which is rewritten as in
(\ref{sommaDeltaYDeltatildeJ2}) with
$\Delta_{j}\tilde{J}_{2}^{(1)}\Delta_{j}\tilde{J}_{2}^{(2)}$ in
place of $\Delta_{j}Y^{(1)}\Delta_{j}\tilde{J}_{2}^{(2)}$, so that
last three terms converge to zero in probability as before. As for
the first term
\begin{equation}\label{sommaDeltatildeJ1DeltatildeJ2}
\Plim
   \sum_{j=1}^{n}\Delta_{j}\tilde{J}_{2}^{(1)}\Delta_{j}\tilde{J}_{2}^{(2)}1_{\{(\Delta_{j}X^{(1)})^{2}\leq
r_{h},(\Delta_{j}\tilde{J}_{2}^{(1)})^{2}\leq
4r_{h}\}}1_{\{(\Delta_{j}X^{(2)})^{2}\leq
r_{h},(\Delta_{j}\tilde{J}_{2}^{(1)})^{2}\leq 4r_{h}\}},
\end{equation}
we remark that it is bounded in absolute value by
$$
\Plim\sum_{j=1}^{n}|\Delta_{j}\tilde{J}_{2}^{(1)}|
1_{\{(\Delta_{j}\tilde{J}_{2}^{(1)})^{2}\leq
4r_{h}\}}|\Delta_{j}\tilde{J}_{2}^{(2)}|1_{\{(\Delta_{j}\tilde{J}_{2}^{(2)})^{2}\leq
4r_{h}\}}
$$
$$
\leq
\Plim\sqrt{S_n^{(1)}} \sqrt{S_n^{(2)}}=0,
$$
by Lemma \ref{SumDXddoveleqrrizero}.

 \vspace{-0.7cm}\qed

\end{document}